\documentclass[11pt]{report} 

\def \dfll {\leaders \hbox to 1em {\hss.\hss}\hfill}
\def\beq {\begin{equation}}
\def\eeq {\end{equation}}

\usepackage{diagrams}

\newarrow{MyInto}{vee}---{vee}
\newarrow{MyOnto}----{doublevee}
\newarrow{MyTo}----{vee}

\newtheorem{theorem}{Theorem}
\newtheorem{propos}[theorem]{Proposition}
\newtheorem{coro}[theorem]{Corollary}
\newtheorem{conjecture}[theorem]{Conjecture}

\newtheorem{qes}{Question}
\newtheorem{lemma}{Lemma}
\def\blm{\begin{lemma}}
\def\elm{\end{lemma}}
\def\bdf{\begin{Defi}}
\def\edf{\end{Defi}}
\def\btm{\begin{theorem}}
\def\etm{\end{theorem}}
\def\bpp{\begin{propos}}
\def\epp{\end{propos}}
\def\bQ {\begin{qes}}
\def\eQ {\end{qes}}
\def\btm{\begin{theorem}}
\def\etm{\end{theorem}}
\def\ben{\begin{enumerate}}
\def\een{\end{enumerate}}

\def\byg{\begin{Young}}
\def\eyg{\end{Young}}
\def\bsyg{\begin{SmallYoung}}
\def\esyg{\end{SmallYoung}}
\def\btid{\begin{Tabloid}}
\def\etid{\end{Tabloid}}
\def\bstid{\begin{SmallTabloid}}
\def\estid{\end{SmallTabloid}}

\def\wave{$\raise+3.3pt\hbox{$\wr$}\mkern-5mu\raise-3.3pt\hbox{$\wr$}$}

\mathchardef\xx="167F

\def\dwa{$\raise+7pt\hbox{$\langle\mkern-3.5mu\langle$}
\mkern-8.5mu\raise-0pt\hbox{$\rangle\mkern-3.5mu\rangle$}
\mkern-13mu\raise-10pt\hbox{$\Downarrow$}$}

\def\sdwa{$\raise+5pt\hbox{$\scriptstyle\rangle\mkern-3.5mu\rangle$}
\mkern-9.55mu\raise-0pt\hbox{$\scriptstyle\langle\mkern-3.5mu\langle$}
\mkern-6.95mu\raise-5pt\hbox{$\scriptstyle\rangle\mkern-3.5mu\rangle$}
\mkern-11.9mu\raise-11.1pt\hbox{$\scriptstyle\bigtriangledown$}
$}

\def\into{\hookrightarrow}
\def\ep{\ \hfill{\rule {2.5mm}{2.5mm}}\smallskip}
\newcommand{\lbl}[1]{\label{#1}}

\def \Cob {{\cal C}ob_3}

\def\mod{\,{\rm mod}\,}

\def\trig{\,\scriptstyle\triangle}
\def\nil#1{\vbox{\ialign{##\crcr\noalign
             {\kern-3pt\nointerlineskip}{$\hfil\,\trig\hfil$} \crcr\noalign
             {\kern-.5pt\nointerlineskip}
             $\hfil\displaystyle{#1}\hfil$\crcr}}\!}

\def\opc{{\scriptstyle =}}

\def\drulefill{$\mathord\opc\mkern-3mu\cleaders
                  \hbox{$\mkern-2mu\mathord\opc\mkern-2mu$}\hfill\mkern-2mu \mathord\opc$}

\def\dov#1{\vbox{\ialign{##\crcr\noalign
             {\kern-3pt\nointerlineskip}\drulefill \crcr\noalign
             {\kern0pt\nointerlineskip}
             $\hfil\displaystyle{#1}\hfil$\crcr}}}

\def\zer#1{\vbox{\ialign{##\crcr\noalign
             {\kern-3pt\nointerlineskip}
             {$\hfil\;\,\scriptscriptstyle\oslash\hfil\!$} \crcr\noalign
             {\kern.5pt\nointerlineskip}
             $\hfil\displaystyle{#1}\hfil$\crcr}}\!}

\def\ste#1{\vbox{\ialign{##\crcr\noalign
             {\kern-3pt\nointerlineskip}
             {$\hfil\;\,\scriptscriptstyle\bowtie\hfil\!$} \crcr\noalign
             {\kern.5pt\nointerlineskip}
             $\hfil\displaystyle{#1}\hfil$\crcr}}\!}

\def\lb{\lz\!\![}
\def\rb{\rz\!\!]}

\newcommand{\emptystuff}[1]{{}}

\newcommand{\TO}[2]{\stackrel {\mbox{#1}}{\hbox to #2pt{\rightarrowfill}}}

\def\thrafill{$\mathsurround=0pt \mathord- \mkern-6mu 
\cleaders\hbox{$\mkern-2mu
\mathord- \mkern-2mu$}\hfill \mkern-6mu\mathord\twoheadrightarrow$}
\newcommand {\onto} [1]{\hbox to #1pt{\thrafill}}

\usepackage{amsfonts,
            amssymb, young}

\input{epsf.sty}

\newcommand{\Z}{{\mathbb Z}}

\newcommand{\R}{{\mathbb R}}

\newcommand{\N}{{\mathbb N}}

\newcommand{\Aa}{{\mathbb A}}

\newcommand{\F}{{\mathbb F}}
\newcommand{\FF}{{\sf F}\mkern -8mu{\sf F}}

\newcommand{\ff}{{\sf f\mkern -4.4mu f}}

\newcommand{\id}{{\mathbb I}}

\def \sp {{\mathfrak s}{\mathfrak p}}
\def \sl {{\mathfrak s}{\mathfrak l}}

\newcommand{\trA}[1]{\!\!\raise-8pt\hbox{\bsyg #1 \cr \esyg}}
\newcommand{\trB}[1]{\,\, \cdots \!\!\!\raise-8pt\hbox{\bsyg #1 \cr \esyg}} 
\newcommand{\trC}[1]{\,\,\raise5pt\hbox{$\cdots$}\!\!\! 
     \raise2pt\hbox{\bsyg #1\cr \esyg}}

\newcommand{\posit}[2]
{\raise -1.4ex\hbox{${\textstyle #1}\atop {\stackrel{\uparrow}{#2}}$}}

\def \lz  {\langle}
\def \rz  {\rangle}

\newcommand{\head}[1]{\medskip \begin{center}{\Large \sc #1}\end{center}}

\newcommand{\ext}[1] {\mbox{\raisebox{.4ex}{$\bigwedge^{\!#1}$}}\mkern-1mu}
\newcommand{\sext}[1] 
{\mbox{\raisebox{.2ex}{$\scriptstyle \bigwedge^{\!#1}$}}\mkern-1mu}

\def\Ii {\mbox{\raise .4 ex\hbox{$\int$}$\!\! I$}}

\begin{document}

\vspace*{0cm}

\begin{center}

\section*{\LARGE Resolutions of $p$-Modular TQFT's and \\
Representations of Symmetric Groups\footnote{
2000 Mathematics Subject Classification: Primary 57R56, 57M27; Secondary 
 17B10, 17B37, 17B50, 20B30, 20C30, 20C20.}}

\bigskip

\bigskip

\bigskip

{\Large Thomas Kerler}\\
\medskip

September 2001
 \vspace*{1.2cm}


\end{center}


\medskip

\begin{center}
{\Large \sc Contents}\
\smallskip

\parbox[t]{11cm}{
1. Introduction and Survey of Results \dotfill\pageref{S1}

2. Frohman-Nicas  TQFT's over $\Z$\dotfill\pageref{S2}

3. Lefschetz Decompositions and Specht Modules\dotfill\pageref{S3}

4. $\F_p=\Z/p\Z$-Reductions and the Sequences  ${\cal C}_{p,k}$\dotfill\pageref{S4}

5. Exactness of   ${\cal C}_{p,k}$ \dotfill\pageref{S5}

6. Characters, Dimensions, and the Alexander Polynomial\dotfill\pageref{S6}

7. Johnson-Morita Extensions \dotfill\pageref{S7}


}
\end{center}

\

\medskip

\

\medskip

\

\head{1. Introduction and Survey of Results}\lbl{S1}

In this article we bring together several areas of representation theory
in a series of interrelated results. The first is the rather established
theory of $p$-modular representations of the symmetric groups, 
 followed by the representation theory of groups of Lie type over the finite 
field $\F_p$, and, finally, the area of topological quantum field theories
(TQFT's) over $\F_p$. The latter have been our original motivation since
they appear as constant order reduction in the Reshetikhin Turaev Theories.
In fact, as we shall outline in more detail at the end of in this section,
the identities we will find in this article will, for example, imply relations 
between the Lescop  invariant and the Reshetikhin Turaev invariant for 
3-manifolds for $p=5$.

Our results in each area concern  resolutions and expansions of 
$p$-modular representations and invariants into their respective 
characteristic zero counterparts, and, thus, naturally build on each other. 
Let us state the results in order, beginning with the case of the symmetric 
groups.

\paragraph{The Symmetric Groups:}
The representation theory of the symmetric group $S_n$ in $n$ letters 
over ${\mathbb Q}$ (or ${\mathbb Z}$) is well known. The theory is
semisimple and the simple representations are isomorphic to the Specht
modules ${\cal S}^{\tau}$, where $\tau$ is a Young diagram with $n$
boxes. They have a natural basis given by Young tablaux, and the $S_n$-action
preserves the free $\Z$-modules (or lattices) ${\cal S}^{\tau}_{\Z}$ 
generated by these basis vectors. Passing to 
${\cal S}^{\tau}_p=\raise .3ex\hbox{${\cal S}^{\tau}_\Z$}\! \big/\!
\raise -.3ex\hbox{$p{\cal S}^{\tau}_\Z$}$ we thus obtain representations
of the same rank over the field  $\F_p$ for any given 
prime $p\geq 3$. The ${\cal S}^{\tau}_p$, however, 
are not longer irreducible, but they have a unique simple quotient
${\cal D}^{\tau}_p$ obtained from canonical inner forms on the Specht modules,
see \cite{Jam78}. The representations ${\cal S}^{\tau}_\Z$ and 
${\cal D}^{\tau}_p$ define characters  $\chi^{\tau}$ and
$\varphi_p^{\tau}$ on $S_n$ with values in $\Z$ and $\F_p$, respectively.    
We also denote by $\chi^{\tau}_p$ the $p$-reduction of $\chi^{\tau}$, which may,
alternatively, be thought  of as the character asscociated to ${\cal S}^{\tau}_p$. 

The relationship between the ``ordinary''  representations 
${\cal S}^{\tau}_{p\, \mbox{\em\scriptsize or}\,\Z}$ 
and ``ordinary'' characters $\chi^{\tau}_{p\, \mbox{\em\scriptsize or}\,\Z}$ 
and their  $p$-modular counterparts ${\cal D}^{\tau}_p$ and 
$\varphi_p^{\tau}$ is, even after decades of research, still 
intensely investigated with many open questions remaining. 
While several algorithms exist for expressing the $\chi^{\tau}_p$
in terms of the $\varphi_p^{\tau}$, fewer results exist for the converse
relations, and even fewer results relate the modules themselves.
The exact modular structure of the ${\cal S}^{\tau}_p$ in terms of a
modular ordering of the
${\cal D}^{\tau}_p$-components has only very recently
been uncovered by Kleshchev and Sheth in \cite{KS99} for the special case
of Young diagrams $\tau$ with two rows. 

The first result of this article may be thought of as the inverse relation of 
the result in \cite{KS99}, and it implies a similarly inverse relation for
the characters.

\begin{theorem}\label{thm-SnReso}
Let $p\geq 3$ be a prime and $\tau=[a,b]$ be the two
row Young diagram with row lengths $a$ and $b$ with $0\leq a-b\leq p-2$. 
Consider  the sequence of 
Young diagrams $\tau_j$ given  by $\tau_{2i}=[a+ip,b-ip]$ and 
$\tau_{2i-1}=[b+ip-1, a-ip+1]$, whenever defined. (That is,
$\tau_0=\tau=[a,b]$, $\tau_1=[b+p-1,a-p+1]$, $\tau_2=[a+p,b-p]$, 
$\tau_3=[b+2p-1,a-2p+1]$,  etc.) 
\begin{enumerate} 
\item There is a resolution of the modular, simple representation
${\cal D}^{\tau}_p$ of $S_n$ in terms  ordinary representations given by
an exact sequence of $S_n$-equivariant maps as follows.
\beq\label{eq-thm-Cseq}
\ldots\;{\longrightarrow}\;
{\cal S}^{\tau_2}_p\;\;{\longrightarrow}\;\;
{\cal S}^{\tau_1}_p\;\;{\longrightarrow}\;\;
{\cal S}^{\tau_0}_p\;\;\longrightarrow\;\;{\cal D}^{\tau}_p\;\;\to\;\;0\;. 
\eeq
\item We have the following expansion of $p$-modular characters
in terms of ordinary characters.
\beq\label{eq-charactersX}
\varphi^{\tau}_p\;\;=\;\;\sum_{i\geq 0} (-1)^{i} \chi^{\tau_{i}}_p\;, 
\eeq 
\end{enumerate} 
\end{theorem}
It is an intriguing fact that the maps in the sequence 
(\ref{eq-thm-Cseq}) are powers of generators of ${\mathfrak s}{\mathfrak l}_2$
acting dually on $(\Z^2)^{\otimes n}$, which, as an $S_n$-module, contains the 
permutation modules $M^{\tau}$. The precise action  is constructed 
in Corollary~\ref{cor-Cseq} of Section~4, where we prove that it yields a
well defined complex.

The proof of exactness of this complex uses the results
of \cite{KS99} and occupies most of Section~5.
A generalization of Theorem~\ref{thm-SnReso} to 
$n$-row diagrams using a dual ${\mathfrak s}{\mathfrak l}_n$ is likely
to yield more insights the structure of $p$-modular Specht modules for general
Young diagrams. 
The character identity (\ref{eq-charactersX}) is stated as an immediate
consequence in Corollary~\ref{cor-characters} of Section~6.

\paragraph{The Symplectic Groups:}
The symmetric groups typically appear as or within Weyl groups of groups
of Lie type. In this article we are particularly interested in
 the symplectic groups ${\rm Sp}(2g,{\mathbb M})$ where 
${\mathbb M}$ may be ${\mathbb R}$, $\Z$ or ${\mathbb F}_p$. 
There are obvious generalizations of our results to most 
other groups of Lie type, which we leave to the reader.  

Let $H=H_1(\Sigma_g,\Z)$ be the fundamental representation of 
${\rm Sp}(2g,\Z)$  with symplectic
basis $\langle a_1,\ldots,a_g,b_1,\ldots,b_g\rangle$, and denote by
$V(\varpi_k)\subset \ext k H$ the subrepresentation generated by the
highest weight vector 
$a_1\wedge\ldots \wedge a_k$. Over ${\mathbb R}$ this is the irreducible 
representation of fundamental
heighest weight $\varpi_j=\epsilon_1+\ldots+\epsilon_j$.
Denote by $V_p(\varpi_j)$ the respective $p$-modular reduction, and 
by $\dov V_p(\varpi_j)$ the unique irreducible subquotient over ${\mathbb F}_p$
generated by the same highest weight vector.
Let $g-p+2\leq l\leq g$ and  $\hat l=2(g-p+1)-l$. 
From Theorem~\ref{thm-SnReso} now we derive resolutions of 
Sp-representations
that are given
by exact sequences as follows. 

\beq\label{eq-SpReso}
\ldots
\;\to\;
V_p(\varpi_{l-4p})
\;\to\;
V_p(\varpi_{{\hat l}-2p})
\;\to\;
V_p(\varpi_{l-2p})
\;\to\;
V_p(\varpi_{\hat l})
\;\to\;
V_p(\varpi_l)
\;\to\;
\dov V_p(\varpi_l)
\;\to\;
0
\,.
\eeq

Evidently, (\ref{eq-SpReso}) implies similar expansions of 
${\rm Sp}$-characters, and such expansions also exist for other
groups of Lie type. The generalization more intersting to us,
which in fact implies (\ref{eq-SpReso}), are resolutions of 
topological quantum field theories (TQFT's).

In Lemmas~\ref{lm-upsilambda} and \ref{lm-induce} of Section~2
we establish the relation between the ${\rm Sp}(2g,\Z)$ weight spaces of 
$\ext *H$ with the $S_k$-module $(\Z^2)^{\otimes k}$. 
In Lemma~\ref{lm-Spe-action} we also give the explicit action of the
Serre generators of ${\mathfrak s}{\mathfrak p}_{2g}$ in the
language of the permutation modules. The relation between the weight 
spaces of $V(\varpi_k)$ and the Specht modules and the respective explicit 
actions of the ${\mathfrak s}{\mathfrak p}_{2g}$-generators on the 
Young diagrams is derived in Section~3. 
 
\medskip

\paragraph{Homological TQFT's:}
Recall that a TQFT is a functor ${\cal V}:\Cob\to {\mathbb M}{\rm -mod}$ from a category of
2+1-dimensional cobordisms into a category of free modules over a commutative ring ${\mathbb M}\,$.
Specifically, $\cal V$ assigns to a surface $\Sigma$ a free ${\mathbb M}\,$-module ${\cal V}(\Sigma)$ 
and to a cobordism between two surfaces an  ${\mathbb M}\,$-linear map between the respective 
${\mathbb M}\,$-modules. 

In \cite{FroNic92} Frohman and Nicas construct a TQFT ${\cal V}^{FN}$ over ${\mathbb M}=\Z$,  
where the free $\Z$-module for a surface is given by the (co)homology of its Jacobian,
specifically,  ${\cal V}^{FN}(\Sigma)=H^*({\rm Hom}(\pi_1(\Sigma),U(1)),\Z)
=\ext *H_1(\Sigma,\Z)\,$.
The basics of this construction are reviewed in the beginning of Section~2. 
Further, in Section~3, we will recall the decomposition  
of this TQFT into its irreducible components ${\cal V}^{(j)}_\Z$,
with $j=1,2,\ldots\,$, of ${\cal V}^{FN}$. They are again TQFT's over $\Z$ and  are obtained 
in \cite{Ker01} from a dual Lefschetz ${\mathfrak s}{\mathfrak l}_2$-action.

Composing ${\cal V}^{(j)}_Z$ with  the canonical (rank-preserving) functor  
$\,{\mathbb Z}{\rm -mod}\,\to\!\!\!\!\to\,{\mathbb F}_p{\rm -mod}\,$ for primes $p\geq 3$
we obtain a family of TQFT's ${\cal V}^{(j)}_p$ over ${\mathbb M}={\mathbb F}_p$. As before, the 
$p$-modular TQFT's are generally highly reducible. However, they have a unique irreducible 
TQFT subquotient $\dov {\cal V}^{(j)}_p$. 

\begin{theorem}\label{thm-TQFTreso} For any prime 
  $p\geq 3$   and integer $k$ with $0<k<p$ we have an exact sequence 
 of natural transformations of TQFT's 
\beq\label{eq-TQFTreso}
\ldots\,\to\,{\cal V}^{((i+1)p+k_{i+1})}_p\,\to\,
{\cal V}^{(ip+k_i)}_p\,\to\,\ldots\,\to\,
{\cal V}^{(2p-k)}_p\,\to\,
{\cal V}^{(k)}_p\,\to\,\dov {\cal V}^{(k)}_p\,\to\,0\;, 
\eeq
where we set  $k_i=k$ for even $i$ and $k_i=p-k$ for odd $i$. 
\end{theorem}

The sequence is constructed from its symmetric group summands in Corollary~\ref{cor-seqTQFT}
of Section~4. Exactness follows in Section~5 from the respective results for 
Specht modules.

In order to see why (\ref{eq-SpReso}) is indeed a special case of 
Theorem~\ref{thm-TQFTreso} observe that
the mapping class group $\Gamma_g$ is identical with the group of 
invertible cobordisms in $\Cob$ from a surface 
$\Sigma_g$ to itself. Hence, any TQFT $\cal V$ entails for every $g$ 
a representation  of $\Gamma_g$ on ${\cal V}(\Sigma)$,
which, in the case of ${\cal V}^{FN}$, factors through the symplectic 
quotient $\Gamma_g\to\!\!\!\!\to {\rm Sp}(2g,\Z)$. The sequence in 
(\ref{eq-SpReso}) is now obtained from Theorem~\ref{thm-TQFTreso} 
by evaluation on a particular surface and 
using  the identifications  of ${\rm Sp}(2g,\Z)$-modules, given by  
${\cal V}^{(j)}_\Z(\Sigma_g)\cong V(\varpi_{g-j+1})$ and 
$\dov {\cal V}^{(j)}_p(\Sigma_g)\cong \dov V_p(\varpi_{g-j+1})\,$.

In the TQFT framework characters are endowed with an interpretion as topological 
invariants of infinite cyclic covers of closed 3-manifolds.  
More precisely, let $\cal V$ be any TQFT, $M$
a closed 3-manifold with $b_1(M)\geq 1$, and
$\varphi:H_1(M)\to\!\!\!\to \Z$ an epimorphism. We define an
 invariant of a pair $(M,\varphi)$ as follows. Pick any two-sided 
surface $\Sigma\subset M$ which is dual to $\varphi\,$, 
and define $C_{\Sigma}=\overline{M-\Sigma}$ seen
as a cobordism from $\Sigma$ to itself. The value
${\cal V}(M,\varphi)=trace({\cal V}(C_{\Sigma}))$ is now independent
of the choice of $\Sigma$. 

As an important example, Frohman and Nicas extracted in \cite{FroNic92} 
the Alexander polynomial $\Delta_{\varphi}(M)$ as a Lefschetz trace 
from ${\cal V}^{FN}$. Translated into the decomposition
of \cite{Ker01} this result says that the invariant ${\cal V}^{(j)}_\Z(M,\varphi)$ 
is the difference of two successive coefficients of the Alexander polynomial. 
Combining this observation with Theorem~\ref{thm-TQFTreso} we derive in the end 
of Section~6 the following relation between these invariants over $\Z$ and the 
respective invariants  $\dov {\cal V}^{(j)}_p(M,\varphi)$
obtained from the {\em irreducible} TQFT's over ${\mathbb F}_p$.

\begin{theorem}\label{thm-pAlexander}  Let $\Delta_{\varphi}(M)\in\Z[x,x^{-1}]$ be the
Alexander Polynomial of a closed, compact, oriented 
3-manifold $M$ with respect to an epimorphism $\varphi:H_1(M)\to\mkern-15mu\to\Z$.

Let $\dov\Delta_{\varphi, p}^{\pm}(M)$ be the reduction of  
$\Delta_{\varphi}(M)$ 
obtained by substituting $x=\pm \zeta_p$, with $\zeta_p$ a $p$-th root of unity,
 and taking the $\Z$ coefficients modulo $p$.
Then 
$$
\dov\Delta_{\varphi, p}^{\mp}(M)\;\;=\;\;\sum_{k=1}^{p-1}(\pm1)^{k-1}[k]_{\zeta_p}
\dov {\cal V}^{(k)}_p(M,\varphi)\;\;\;\;\in\;\;\F_p[\zeta_p]\;. 
$$
\end{theorem}

\paragraph{Johnson-Morita-Extensions and TQFT-Overview:}
Note, that all TQFT's  constructed up to this point have representations of $\Gamma_g$
that factor through ${\rm Sp}(2g,\Z)$,
that is, they contain the Torelli  group ${\cal I}_g$ in their kernels. In Section~7
of this article we will also consider extensions of these TQFT's with slightly 
smaller kernels, at least at the level of representations of $\Gamma_g$.

In  \cite{Mor93} Morita constructs a homomorphism 
$\tilde k: \Gamma_g\to {\rm Sp}(2g,\Z)\ltimes \frac 12 U$ which has a smaller kernel
${\cal K}_g\subset{\cal I}_g$, and whose image ${\bf Q}_g\cong\Gamma_g/{\cal K}_g$ has finite index
in ${\rm Sp}(2g,\Z)\ltimes \frac 12 U$. It extends the Johnson 
homomorphism ${\cal T}_g\to\!\!\!\to U$ previously
constructed in \cite{Joh80}, where $U$ denote the free abelian
group $U\cong\,\raise .4ex\hbox{$\ext 3 H$}\!\Big / \!\raise -.4ex\hbox{$\omega\wedge H$}
\,\cong \,V(\varpi_3)$. As before we denote 
$H=H_1(\Sigma_g)$, considered as an Sp-module, and we let $\omega\in\ext 2H$
be  the standard Sp-invaraiant symplectic form. 

Non-trivial representations of 
 ${\bf Q}_g$ are readily obtained in Proposition~\ref{propos-Umods} of 
Section~7 
as the extension of a pair of representations $V(\varpi_l)$
and $V(\varpi_{l+3})$ using the (up to scale unique) Sp-equivariant map
$U\,\to\,{\rm Hom}(V(\varpi_l), V(\varpi_{l+3}))$. We prove that these
extension also exist for the irreducible $p$-modular representations 
$\dov V_p(\varpi_{l})$. In TQFT language this translates to the following result.

\begin{theorem}\label{thm-JMext} For $0<k<p-3$ there are indecomposable
 representations $\dov {\cal U}^{(k)}_p(\Sigma_g)$ 
of the mapping class group $\Gamma_g$ that factor through 
${\bf Q}_g=\Gamma_g/{\cal K}_g$  but represent 
 the Torelli group ${\cal I}_g$ non-trivially. There is a
short, non-split exact sequence of $\Gamma_g$-equivariant maps as follows: 
\beq\label{eq-JMext}
0\;\to\;
\dov {\cal V}^{(k+3)}_p(\Sigma_g)\;\longrightarrow\;
\dov {\cal U}^{(k)}_p(\Sigma_g)\;\longrightarrow\;
\dov {\cal V}^{(k)}_p(\Sigma_g)\;\to\;0\;\;.
\eeq
\end{theorem}
We know  from Theorem~\ref{thm-KerFib} below 
that $\dov {\cal U}^{(1)}_5$ does on fact extend to a TQFT. 
However, the question whether or how the $\dov {\cal U}^{(k)}_p$ and ${\cal U}^{(k)}_{\Z}$
descend from TQFT's for general $p$ and $k$ is still open, and will be discussed in
future work.

We summarize the different TQFT's of this paper, the constructions connecting them, 
and the included modules of the symmetric groups $S_n$ in the following table:
\bigskip

\begin{tabular}{rclcr}
&${\cal V}^{FN}_\Z$ & Jacobian TQFT \cite{FroNic92}  & $\hookleftarrow$ & Permutation Modules\\
&&                    Fully reducible over $\Z$ \vspace*{-.4cm}&&\\
$SL(2,\R)$-Lefschetz &\sdwa\vspace*{-.3cm}  & & & \\
Decomposition  & & & \\
&${\cal V}^{(k)}_\Z$ & Lefschetz Component, & $\hookleftarrow$ & Specht Modules  ${\cal S}_{\Z}^{\tau}$\\
&&                    irreducible  over $\Z$ \vspace*{-.2cm}&&\\
$p$-Reduction &\sdwa\vspace*{.1cm}  & & & \\
&${\cal V}^{(k)}_p$ & Reducible over $\F_p$ & $\hookleftarrow$ & $p$-Specht Modules  ${\cal S}_p^{\tau}$\\
&& with inner form.   \vspace*{-.4cm}&&\\
Null space quotient &\sdwa\vspace*{.1cm}  & & & \\
&$\dov {\cal V}^{(k)}_p$ & Irreducible over $\F_p$. & $\hookleftarrow$ & Simple $S_n$-modules  ${\cal D}_p^{\tau}$\\
&&  Not $\Z$-liftable, but \vspace*{-.1cm} &&\\
&& resolutions in ${\cal V}^{(j)}_p$'s   \vspace*{-.6cm}&&\\
Johnson-Morita &\sdwa\vspace*{-.3cm}  & & & \\
Extension & & & \\
&$\dov {\cal U}^{(k)}_p$ & Indecomposable with &   &  \\
&&        two factors $\dov {\cal V}^{(k)}_p$ and $\dov {\cal V}^{(k+3)}_p$ \hspace*{-1cm}&&\\
\end{tabular}
\bigskip

\paragraph{Relations to the Reshetikhin Turaev Theory:}  
The original motivation of this article comes from the study of the TQFT's 
${\cal V}_{\zeta_p}^{RT}$ 
constructed by Reshetikhin and Turaev in \cite{ResTur91} from
the quantum group $U_{\zeta_p}({\mathfrak s}{\mathfrak o}_3)$ for a $p$-th root of
unity $\zeta_p$.  Thus, in order to put the TQFT's of this article into their
broader context and illustrate their relevance let us 
 briefly sketch the pertinent relations and 
results from other work in quantum topology.

It is shown by Gilmer in \cite{Gil01}
that, for a restricted set of cobordisms, ${\cal V}_{\zeta_p}^{RT}$ 
can be written as a TQFT over  the
ring of cyclotomic integers $\Z[\zeta_p]$. For the Reshetikhin Turaev invariants of 
closed 3-manifolds  this integrality property was previously proved
in \cite{Mur95, MasRob97}.

Applying the ring reduction $\Z[\zeta_p]\to\!\!\!\to {\mathbb F}_p$, with $\zeta_p\mapsto 1$,
we obtain a TQFT ${\cal V}_p^{I}$ over the finite field ${\mathbb F}_p$ from
${\cal V}_{\zeta_p}^{RT}$ for a given $p$. It is now natural to ask whether or not
there exists a relationship 
between the ${\cal V}_p^{I}$ and the $\dov{\cal V}^{(j)}_p\,$, since they are both
TQFT's over ${\mathbb F}_p$ and are conjectured to share a list of other 
features. 

Identifications between the TQFT's ${\cal V}_p^{I}$ obtained from quantm 
groups on the one hand and the  ones obtained from the homological theory 
${\cal V}_p^{FN}$ (and their Johnson Morita extensions)  
on the other hand will entail many insights into the
relation between classical and quantum invariants as well as establish natural bases
for the quantum theory in the language of the tensor calculus of the symplectic groups. 
Moreover, we expect this to lead to a
deeper understanding of the geometric interpretations of  the higher order terms in the 
cyclotomic integer expansions of the Reshetikhin Turaev theories, which include the
Ohtsuki and, particularly, the Casson-Walker-Lescop invariants.

Since the ${\cal V}_{\zeta_p}^{RT}$ are ``unitarizable'' theories and the 
dimensions of the vector spaces do not match any combination of ordinary Sp-representations 
we must expect the simple quotients $\dov{\cal V}^{(j)}_p$ to enter the relations
rather than the reducible TQFT's ${\cal V}_p^{(j)}$. This is precisely the point where the 
main results of this article are crucially needed. The resolutions of $p$-modular TQFT's provide 
the missing link between the classical invariants defined over $\Z$ and the quantum invariants 
defined over ${\mathbb F}_p$. 

Thus far we are able to understand these relations and their applications in a rather detailed
manner for $p=5$. We state next the main result from  \cite{KerFib}. 

\begin{theorem}[\cite{KerFib}]\label{thm-KerFib} We have the following isormophism of 
$\Gamma_g$-representations over $\F_5$:
\beq\label{eq-TQFTID}
{\cal V}^{I}_{5}(\Sigma_g)\;\;\cong\;\;\dov {\cal U}^{(1)}_5(\Sigma_g)\;\;.
\eeq
In particular, for pairs $(M,\varphi)$ as in Theorem~\ref{thm-pAlexander} we
have 
\beq\label{eq-Lescop}
{\cal V}^{I}_{5}(M,\varphi)\;\;=\;\;\dov {\cal V}_5^{(1)}(M,\varphi)\,+\,
\dov {\cal V}_5^{(4)}(M,\varphi)\;\;=\;\;-2\cdot \lambda_{Lescop}(M)\qquad{\rm mod}\,\,5
\eeq
\end{theorem}

The first equality in (\ref{eq-Lescop}) follows directly from (\ref{eq-TQFTID}). The second
equality is now a consequence of the identification with the Alexander polynomial at a root
of unity from  Theorem~\ref{thm-pAlexander}. It is used in \cite{KerKyoto} to identify 
${\cal V}^I_5(M,\varphi)$ with the ${\mathbb F}_5$-reduction of the Lescop invariant
$\lambda_{Lescop}(M)\,$. This identity is special to $p=5$ and does not hold for $p>5$. 

Observe also that (\ref{eq-Lescop}) implies that the invariant ${\cal V}^I_5(M,\varphi)$ 
is really independent of $\varphi$. This is not surprising since it is also equal to the
lowest oder non-vanishing coefficient of the cyclotomic integer expansion of  
${\cal V}^{RT}_{\zeta_5}(M)$. Nevertheless, this raises an interesting question, namely, which
linear or polynomial combinations of the $\dov {\cal V}_p^{(j)}(M,\varphi)$ are
independent of the choice of $\varphi$ and, thus, yield true invariants of closed 3-manifolds.

\paragraph{Dimensions, Combinatorics, and Asymptotics:}
An important special case for any character formula is the implied  relation for the dimensions
of modules, that is, the characters evaluated at the unit element. The dimensions of the
irreducible modules over ${\mathbb F}_5$ in this article are Fibonacci numbers, while the
dimensions  of the corresponding modules over $\Z$ are Catalan numbers. As a 
result we obtain in (\ref{eq-fibcatideven}) and   (\ref{eq-fibcatidodd})
of Section~6 identities that express the Fibonacci 
numbers as 5-periodic, alternating sums in these Catalan numbers.
Despite the simplicity of these relations they appear to be unknown in the 
available literature. 
  
We also determine in Section~6  the asymptotics of the dimensions of the TQFT modules
for $g\to\infty$, which is summarized in the following lemma.
\begin{lemma}\label{lm-asym} For a fixed prime $p\geq 3$ and fixed $j$, the dimensions 
of the vector spaces grow for large $g$ as follows:
\beq\label{eq-PFasy}
{\rm dim}({\cal V}^{I}_p(\Sigma_g))\;\sim\;{\rm const.}\cdot\|\FF_p\|^g\qquad\;\;{\rm and }
\qquad\;\;{\rm dim}(\dov {\cal V}^{(j)}_p(\Sigma_g))\;\sim\;{\rm const.}\cdot\|\ff_p\|^g
\eeq
where \qquad\qquad
$\displaystyle
\| \FF_p\|=\frac p{4\sin^2(\frac \pi p)}\qquad{\rm and }\qquad \| \ff_p\|=4\cos^2(\frac \pi{2p})\;.
$

Moreover, there are polynomials $R_j(f)\in\Z[f]$ of degree $\,{\rm deg}(R_j)=\frac{p-3}2\,$ with
\beq\label{eq-FfPoly}
\|\FF_p\|\;\;=\;\;R_p(\|\ff_p\|)\;\;.
\eeq
\end{lemma}
The polynomial dependence between the asymptotic behaviors from (\ref{eq-FfPoly}) 
suggests a similar ``polynomial'' dependence between the TQFT's. If the products of
the dimensions are replaced by tensor products of the respective TQFT's, this suggests
that $\sim R_p(\|\ff_p\|)^g$ describes the asymptotics of parts of an $\frac {p-3}2$-fold
tensor product of the $\dov {\cal V}^{(j)}_p$. We are thus led to the  following 
conjecture on the constant order structure of the Reshetikhin Turaev theory, which
 can be verified for the genus=1 mapping class group (${\rm SL}(2,\Z)$) representations.
\begin{conjecture}\label{conj-sums}\ 

The irreducible components of ${\cal V}^I_p$ can be found as irreducible 
components of $\displaystyle \bigotimes^{\frac {p-3}2}\dov {\cal V}^{FN}_p$.
\end{conjecture}

\paragraph{Acknowledgements: } I am much indebted to Alexander Kleshchev for pointing 
out his results in \cite{KS99},  and to Gordon James and Alain Reuter for checking the 
arguments in Section~5. I thank Ronald Solomon, Stephen Rallis, Adalbert Kerber, and 
Charles Frohman for comments and interest. Also, I have benefited much from discussions
with Pat Gilmer on related integral TQFT's and their applications. 
Finally, I would like to express my thanks to Tomotada Ohtsuki and Hitoshi Murakami
for hospitality at the R.I.M.S., Kyoto,  where final parts of this paper were completed,  
and for the opportunity to present and discuss the results.

\head{2. Frohman-Nicas  TQFT's over $\Z$}\lbl{S2} 
\nopagebreak

In \cite{FroNic92} Frohman and Nicas construct a TQFT  that 
allows the interpretation of the Alexander Polynomial as a  weighted 
Lefschetz trace. The vector spaces are the cohomology rings of the
Jacobians of the surfaces. In \cite{Ker01} we extract a natural Lefschetz 
$SL(2,\R)$-action on these spaces, with respect to which this TQFT is 
equivariant.  Let us describe in this sections the basic ingredients of the 
Frohman Nicas theory, reviewing also the conventions of \cite{Ker01}. 

 We denote by 
$\Cob^{22fr}$ the category of {\em evenly framed} cobordisms
between connected standard surfaces $\Sigma_g$.
These are essentially the cobordisms obtained by an even number 
of surgeries, see Lemma~10, \cite{Ker01}. The Frohman Nicas TQFT
is then given as a functor
\begin{equation}\label{eq-FNfunct}
{\cal V}^{FN}\;:
\;\;\Cob^{22fr}\;\longrightarrow\;SL(2,\R)-mod_{\R}\;,
\end{equation}
for which
$$
{\cal V}^{FN}(\Sigma_g)\;\;=\;\;\ext * H_1(\Sigma_g)\;\;
$$
for the surface $\Sigma_g$ of genus $g$. For each surface we 
pick a symplectic basis 
$\{a_1,\ldots, a_g, b_1,\ldots,b_g\}$ for $H_1(\Sigma_g)$ and 
introduce the inner form $\lz \_,\_\rz$, for which this basis
orthonormal. We denote $i_*:H_1(\Sigma_g)\to H_1(\Sigma_{g+1})$
the inclusion map compatible with the choice of bases and the 
construction of $\Sigma_{g+1}$ from $\Sigma_g$ by handle addition. 

The bases  naturally provide lattice bases for $\ext * H_1(\Sigma_g,\Z)$,
given by monomials in the $a_i$ and $b_i$. Moreover, the inner
form extends to to the exterior product to make the monomial
basis an orthonormal one. We also denote by
$\omega_g=\sum_{i=1}^ga_i\wedge b_i\in \ext 2 H_1(\Sigma_g,\Z)$
the symplectic form. 
The $\sl (2,\R)$-action on $\ext * H_1(\Sigma_g)$ 
is given in terms of the standard generators $E$, $F$, and $H$ as 
\beq\lbl{eq-defsl2act}
H\alpha =(j-g)\alpha\quad\; \forall \alpha\in \ext jH_1(\Sigma_g)\;,
\qquad\quad E\alpha =\alpha\wedge\omega_g\;,\qquad\quad 
F = E^*\;. 
\eeq
From the explicit actions of the generators in \cite{Ker01}
we see that the lattices $\ext * H_1(\Sigma_g, \Z)$ are preserved by the maps 
${\cal V}^{FN}(M)$ for cobordisms $M$ in $\Cob^{22fr}$.
Moreover, also the standard generators of $\sl (2,\R)$ 
preserve the lattices. Hence also the universal enveloping algebra  over $\Z$ 
generated by  the operators in (\ref{eq-defsl2act}), which we shall
(abusively) denote by $\sl (2,\Z)$.  We denote  the respective functor into 
$\Z$-modules by:
\begin{equation}\label{eq-FNSL2funct}
{\cal V}_{\Z}\;:
\;\;\Cob^{22fr}\;\longrightarrow\;\sl(2,\Z)-mod_{\Z }\;,
\end{equation}

More specifically,
the extension of the mapping class group in $\Cob^{22fr}$ factors under
${\cal V}_{\Z}$ through a  split $\F_2=\Z/2\Z$ extension of ${\rm Sp}(2g,\Z)$,
acting in the natural way on $\ext * H_1(\Sigma_g, \Z)$. 
We also denote by $\sp(2g,\Z)$ the algebra over $\Z$ generated by 
the standard ${\rm Sp}(2g,\R)$ Lie algebra generators. It is not hard to
show that in this representation it coincides with the group
algebra $\Z[{\rm Sp}(2g,\Z)]$. Besides the 
mapping class groups the other generators of $\Cob^{22fr}$ are the 
handle attachment cobordisms $H^+_g:\Sigma_g\to\Sigma_{g+1}$ and  
$H^-_g:\Sigma_{g+1}\to\Sigma_g$ with actions given by 
\begin{equation}\label{eq-handle}
{\cal V}_{\Z}(H^+_g)\;:\;\;\alpha\,\mapsto\,i_*(\alpha)\wedge a_{g+1}
\qquad\qquad\mbox{and}\qquad  {\cal V}_{\Z}(H^-_g)={\cal V}_{\Z}(H^+_g)^*\;,
\end{equation}
where we use the notation  $i_*$ also for the  inclusion 
$H_1(\Sigma_g)\into H_1(\Sigma_{g+1})$ extended to the exterior powers.  
\medskip

The structure of ${\cal V}_{\Z}$ and the dual $\sl (2,\Z)$-action
can be better understood if we decompose the lattices ${\cal V}_{\Z}(\Sigma_g)$
according to $\sp(2g,\Z)$-weights. We write an 
$\sp (2g)$-weight in the standard basis 
$\lambda=\sum_{i=1}^g\lambda_i\epsilon_i$ as given in 
 \cite{GooWal98}. Specifically, we have $ha_i=\lz \epsilon_i, h\rz a_i$
and $hb_i=-\lz \epsilon_i, h\rz b_i$ for $h\in {\mathfrak h}$, the 
diagonal matrices in $\sp(2g,\Z)$.

 We denote by $\nabla_g=\{\lambda\in {\mathfrak h}^*:\;\lambda_i\in\{-1,0,1\}\}$
the set of possible weights of vectors in $\ext * H_1(\Sigma_g)$. This yields
a decomposition into weight spaces denoted as follows.
\begin{equation}
{\cal V}_{\Z}(\Sigma_g)\;\;=\;\;\bigoplus_{\lambda\in\nabla_g}
{\cal W}_{\Z}(\lambda,g)
\end{equation}
For a given weight $\lambda\in \nabla_g$ let 
$N({\lambda})=\{ i:\,\lambda_i=0\}\subset\{1,\ldots,g\}$ and 
$n(\lambda)=|N({\lambda})|$.
A special vector $w(\lambda)\in {\cal W}_{\Z}(\lambda,g)$ is given by 
$$
w(\lambda)=w_1(\lambda_1)\wedge\ldots\wedge w_g(\lambda_g)
\;\in\ext {g-n(\lambda)} H_1(\Sigma_g)
\qquad\quad\mbox{where}\quad \left\{
\begin{array}{rcl}w_i(1)&=&a_i \\ w_i(-1)&=&b_i\\ w_i(0)&=&1
\end{array}\right. \;\;. 
$$
Let $e_{\pm}$ be generators
of a 2-dimensional lattice $\lz e_-, e_+\rz_{\Z}$, and  
$L_{\Z}^n=\lz e_-, e_+\rz_{\Z}^{\otimes n}$ the lattice of rank $2^n$. Given 
$N(\lambda)=\{j_1,\ldots, j_{n(\lambda)}\}\,$ with
$j_1<\ldots < j_{n(\lambda)}$ we thus define maps
\begin{eqnarray}\label{eq-upsilon}
\Upsilon_{\lambda}\,&:
&\;\;\;\quad  L_{\Z}^{n(\lambda)}\,
\;\;\stackrel{\cong}{-\!\!\!-\!\!\!-\!\!\!\longrightarrow}
\;\;\, {\cal W}_{\Z}(\lambda,g)\\
\;&:&\;\;e_{\epsilon_1}\otimes\ldots\otimes e_{\epsilon_{n}}
\;\mapsto\;w(\lambda)\wedge o_{j_1}(\epsilon_1)\wedge\ldots\wedge 
o_{j_{n}}(\epsilon_{n})
\nonumber
\\
&&\mbox{where}\quad
n=n(\lambda)\quad\mbox{and}
\quad\quad\left\{\begin{array}{l}
o_j(+)=a_j\wedge b_j\\
o_j(-)= 1\nonumber
\end{array}
\right.\;\;\;. 
\end{eqnarray}
The lattices $L_{\Z}^n$ also have a natural inner product for which
the monomials $e_{\epsilon_1}\otimes\ldots\otimes e_{\epsilon_n}$  form
an orthonormal basis, and carry a natural $\sl (2,\Z)$-action 
given by $Ee_-=e_+$, $Fe_+=e_-$, $He_{\pm}=\pm e_\pm$, and 
$Ee_+=Fe_-=0$, for which also $E^*=F$ and $H^*=H$. 
\begin{lemma}\label{lm-upsilambda}
The $\Upsilon_{\lambda}$ are $\sl (2,\Z)$-equivariant isomorphisms
of lattices with inner forms. 
\end{lemma}

{\em Proof:} They are obviously isomorphisms of lattices as they map 
orthonormal bases to each other. It is also easy to see that the 
$H$-weight for a monomial is $-n(\lambda)+2\sum_{i=1}^{n(\lambda)}\epsilon_i$
on both sides of (\ref{eq-upsilon}). Now, $E$ is multiplication with
$\omega_g=\sum_ja_j\wedge b_j$. Clearly, $w(\lambda)\wedge a_j\wedge b_j=0$
if $j\not\in N(\lambda)$ so we multiply only with 
$\omega_{\lambda}=\sum_{i=1}^{n(\lambda)}a_{j_i}\wedge b_{j_i}=
\sum_{i=1}^{n(\lambda)}o_{j_i}(+)$. The $E$-equivariance then follows from
$o_j(+)\wedge o_j(+)=0$ and $o_j(-)\wedge o_j(+)=o_j(+)$. $F$-equivariance
follows from $F^*=E$.
\ep
\smallskip

Although not of immediate necessity 
for the main result let us record here also how the
actions of morphisms ${\cal V}_{\Z}(M)$ on the lattices $L_{\Z}^n$ look
like. We give them in terms of generators  of $\sp(2g,\Z)$ and the
handle attaching maps. 
We introduce  $\sl (2,\Z)$-equivariant morphisms
\begin{equation}
coev_k\;=\;\id_{k-1}\otimes coev\otimes \id_{n-k}\;\;:\;\;
\qquad L_{\Z}^n\;\;\longrightarrow \;\;
L_{\Z}^{n+2}
\qquad\qquad\qquad\qquad
\end{equation}
$$
\mbox{where}
\qquad coev\,:\; 1\,\mapsto\, (e_-\otimes e_+-e_+\otimes e_-)
\qquad\qquad 
\mbox{and}\quad \id_n\,=\,id_{L^n_{\Z}}\;.
$$
$$
\mbox{and}\qquad ev_k= - (coev_k)^*:\;\;L_{\Z}^{n+2}\longrightarrow L_{\Z}^n
\qquad\qquad\qquad
$$
The obey relation
$$
ev_{k\pm 1}\circ coev_k\;=\;\id \qquad\qquad\mbox{and}\qquad\qquad
ev_{k}\circ coev_k\;=\; -2\cdot\id\;. 
$$
For a  map $\phi:{\cal W}_{\Z}(\lambda,g)\to {\cal W}_{\Z}(\lambda',g')$ 
we denote $\phi^{\Upsilon}=\Upsilon_{\lambda'}^{-1}\circ \phi\circ
\Upsilon_{\lambda}: L_{\Z}^{n(\lambda)}\to L_{\Z}^{n(\lambda')}\,$. 
Moreover, as in \cite{GooWal98}, we denote the standard standard
generators $e_{\alpha_i}$ and $f_{\alpha_j}$ of $\sp (2g,\Z)$ for simple roots 
$\alpha_i=\epsilon_i-\epsilon_{i+1}$ for $i<g$ and $\alpha_g=2\epsilon_g$.
We have $e_{\alpha_i}a_{i+1}=a_i$, $e_{\alpha_i}b_i=-b_{i+1}$, and 
$e_{\alpha_i}v=0$ for all other basis vectors $v$ if $i<g$. Also
$e_{\alpha_g}b_g=e_g$ and $e_{\alpha_g}v=0$ for all others. This
further determines the action of the other generators with 
$f_{\alpha_i}=e_{\alpha_i}^*$. 

Denote now by $e_{\alpha,\lambda}$ the restriction 
$e_{\alpha}:{\cal W}_{\Z}(\lambda,g)\to {\cal W}_{\Z}(\lambda+\alpha,g)$,
where we put $e_{\alpha,\lambda}=0$ of $\lambda+\alpha\not \in \nabla_g$. 
Let us also denote the restriction of the handle attaching map
 $H_{\lambda}={\cal V}_{\Z}(H^+_g):
{\cal W}_{\Z}(\lambda,g)\to {\cal W}_{\Z}(\lambda+\epsilon_{g+1},g+1)$. 
The following is the result of a straight forward calculation. 
\begin{lemma}\label{lm-Spe-action} 
For $\lambda\in\nabla_g$ and $n=n(\lambda)$ we have, when $i<g$, 
$$
e_{\alpha_i,\lambda}^{\Upsilon}\;\;=\;\;
\left\{
\begin{array}{cll}
\quad\id_n\quad & \mbox{for\ \ \ } (\lambda_i,\lambda_{i+1})=(0,1)&\\
-\id_n&  \mbox{for\ \ \ } (\lambda_i,\lambda_{i+1})=(-1,0)&\\
coev_k& \mbox{for\ \ \ } (\lambda_i,\lambda_{i+1})=(-1,1)&
        \mbox{with\ \ \ } j_{k-1}<i<j_k\\
-ev_k& \mbox{for\ \ \ } (\lambda_i,\lambda_{i+1})=(0,0)&
        \mbox{with\ \ \ } i=j_k,\;i+1=j_{k+1}\\
0& \mbox{elsewise}&\\
\end{array}\right. \;\;, 
$$ 
$$
e_{\alpha_g,\lambda}^{\Upsilon}\;\;=\;\;
\left\{
\begin{array}{cll}
\quad\id_n\quad & \mbox{for\ \ \ } \lambda_g\,=\,-1&\\
0& \mbox{elsewise}&\\
\end{array}\right.\;\;, 
\qquad\mbox{and}\qquad H^{\Upsilon}_{\lambda}=\id_n
\qquad.\qquad\qquad\qquad
$$
\end{lemma}

Another prominent, and for our purposes more important, action on
the lattices ${\cal V}_{\Z}(\Sigma_g)$  is that
of subgroups of the Weyl group ${\mathfrak W}_g\cong  (\F_2)^g\rtimes S_g$
of ${\rm Sp}(2g,\Z)$, where $S_g$ denotes the symmetric group in $g$ letters.
The $j$-th  $\F_2$-generator of ${\mathfrak W}_g$ acts on weights by
changing  the sign of $\lambda_j$. It is realized as a subgroup 
$\widehat{\mathfrak W}_g\subset {\rm Sp}(2g,\Z)$ by an extension
$1\to {\mathbb F}_2^g\to\widehat{\mathfrak W}_g\to {\mathfrak W}_g \to 1$, with 
$\widehat{\mathfrak W}_g\cong  (\F_4)^g\rtimes S_g$. 
The $\F_4$-generators are given 
by the ``$S$-matrices'' ${\sf  S}_j\in {\rm Sp}(2g,\Z)$, see \cite{Ker01}, defined by
${\sf S}_ja_j=-b_j$, ${\sf S}_jb_j=a_j$  and ${\sf S}_ja_i=a_i$ and 
${\sf S}_jb_i=b_i$ for $i\neq j$. 
We specify two relevant representations of $\widehat{\mathfrak W}_k$:
\begin{enumerate}
\item $L^k_{\Z}\cong (\Z^2)^{\otimes k}$: This action factors through 
$\widehat{\mathfrak W}_k\to S_k$ the symmetric group, which acts
canonically on the lattice by permutation of factors. 
\item $M^k_{\Z}\cong (\Z^2)^{\otimes k}$: Here the $j$-th $\F_4$-factor is
represented by the matrix ${\sf S}_j=\scriptscriptstyle 
\left[\matrix{0 & 1 \cr -1 & 0}\right ]$ acting on the $j$-th factor
of the tensor product. The action of ${\sf S}_k$ on $M^k_{\Z}$ is the 
canonical representation multiplied by the alternating representation,
i.e., $\sigma(v_1\otimes\ldots \otimes v_k)=sign(\sigma)
v_{\sigma^{-1}(1)}\otimes\ldots \otimes v_{\sigma^{-1}(k)}$.
\end{enumerate} 
We have natural subgroups 
$\widehat{\mathfrak W}_{g-n}\times \widehat{\mathfrak W}_n\subset 
\widehat{\mathfrak W}_g$, for which the right coset 
$$
C^g_n\;=\; \raise3pt\hbox{$\widehat{\mathfrak W}_g$}\Big / 
\raise-3pt\hbox{$\widehat{\mathfrak W}_{g-n}\times \widehat{\mathfrak W}_n$}
\;=\; \raise3pt\hbox{$S_g$}\Big / 
\raise-3pt\hbox{$S_{g-n}\times S_n$}\;\;,
$$    
is identified with the set of subsets $A\subset \{1,\ldots, g\}$ of size $|A|=n$.  
We denote
\begin{eqnarray}\label{eq-Wnotation}
{\cal W}_{\Z}(n,g)&\,=\,&\bigoplus_{A\in C^g_n}  {\cal W}_{\Z}(A,g) 
\,=\,\bigoplus_{\lambda\in\nabla_g: n(\lambda)=n}  {\cal W}_{\Z}(\lambda,g) \;, 
\qquad\qquad\qquad
\\
\ \ &\mbox{where}&\qquad
{\cal W}_{\Z}(A,g)\,=\!\!\bigoplus_{\lambda\in\nabla_g:N(\lambda)=A}  {\cal W}_{\Z}(\lambda,g) \;. 
\end{eqnarray}
Clearly, the summands of ${\cal V}_{\Z}=\bigoplus_n{\cal W}_{\Z}(n,g)$ are 
invariant under the $\widehat{\mathfrak W}_g$-action for each $n$. These
subrepresentations are identified next as induced representations.
\begin{lemma}\label{lm-induce} 
For every $n$ with $0\leq n\leq g$ there 
is a natural isomorphism of $\widehat{\mathfrak W}_g$-modules
$$
\Upsilon\;:\;\;
Ind^{\,\,\widehat{\mathfrak W}_g}_{\widehat{\mathfrak W}_{g-n}\times \widehat{\mathfrak W}_n}
\Bigl(M^{g-n}_{\Z}\otimes L^n_{\Z}\Bigr)
\;\stackrel{\cong}{\longrightarrow}\;\;{\cal W}_{\Z}(n,g)\;. 
$$
This map is an $\sl(2,\Z)$-equivariant isometry. 
\end{lemma}

{\em Proof:} For $N_n=\{g-n+1,\ldots,g\}$ the 
$\widehat{\mathfrak W}_{g-n}\times \widehat{\mathfrak W}_n$-module is
readily identified via the isomorphism 
$$
M^{g-n}_{\Z}\otimes L^n_{\Z}\;\stackrel{\cong}{\longrightarrow}
\;{\cal W}_{\Z}(N_n,g)\;:\;\;
e_{\lambda_1}\otimes\ldots\otimes e_{\lambda_{n-g}}\otimes l
\;\mapsto\;\Upsilon_{[\lambda_1,\ldots,\lambda_{n-g},0,\ldots,0]}(l)\;.
$$
with the submodule ${\cal W}_{\Z}(N_n,g)\subset\ext * H_1(\Sigma_g)$, where
the action is defined by restricting the action of $\widehat{\mathfrak W}_g$.
We next define a natural section 
\beq\label{eq-section}
\pi\;:\; \; C^g_n\;\longrightarrow\; S_g\subset \widehat{\mathfrak W}_g\;\;:
\;\;\;A\;\mapsto\;\pi_A
\eeq
as follows.  Let  
$A\subseteq\{1,\ldots,g\}$ with $n=|A|$. There is a 
unique permutation $\pi_A\in S_g$ such that
$A=\pi_A(\{g-n+1,\ldots,g\})$, $\pi_A(1)<\pi_A(2)<\ldots<\pi_A(g-n)$,
and $\pi_A(g-n+1)<\pi_A(g-n+2)<\ldots<\pi_A(g)$. 
Clearly, we have
$$
\pi_A:{\cal W}_{\Z}(N_n,g)\;\stackrel{\cong}{\longrightarrow}\;
{\cal W}_{\Z}(A,g)\;. 
$$
The induced representation by definition the space of all maps
$f: \widehat{\mathfrak W}_g \to M^{g-n}_{\Z}\otimes L^n_{\Z}$ such that
$f(\sigma \eta)=\eta^{-1}.(f(\sigma))$
for $\eta\in \widehat{\mathfrak W}_{g-n}\times \widehat{\mathfrak W}_n$,
equipped with the left regular 
$\widehat{\mathfrak W}_g$-action $(\sigma f)(\sigma')=f(\sigma^{-1}\sigma')$. Now, every
$\sigma\in \widehat{\mathfrak W}_g$ has a unique decomposition 
$\sigma=\pi_{\sigma(N_n)}\eta_\sigma$, with 
$\eta_\sigma\in \widehat{\mathfrak W}_{g-n}\times \widehat{\mathfrak W}_n$.
Thus we may identify the induced representation with the space of
maps $\overline f: C_n^g\to M^{g-n}_{\Z}\otimes L^n_{\Z}$, setting
$\overline f (A)=f(\pi_A)$  and hence $f(\sigma)=\eta_\sigma^{-1}.\overline f(\sigma(N_n))$. 
The isomorphism is now given by 
\beq\label{eq-isomind}
Ind\,\cong\,Map(C_n^g,\, M^{g-n}_{\Z}\otimes L^n_{\Z})\;\longrightarrow\;{\cal W}_{\Z}(n,g)\;:\;\;
\overline f \;\;\mapsto\;\;\Upsilon(\overline f)=
\bigoplus_{A\in C^g_n}\pi_A(\overline f(A))\;.
\eeq
An inverse is obtained by mapping $v\in {\cal W}_{\Z}(A,g)$ to
$\pi^{-1}_A(v)\otimes \delta_{\{A\}}$, where  $\delta_{\{A\}}(B)=1$
for $A=B$ and 0 elsewise. In order to show that it is equivariant let
$\overline f$ be  an arbitrary map $C^g_n\to {\cal W}_{\Z}(N_n,g)$, 
$\sigma\in \widehat{\mathfrak W}_g$, and $ \overline f'=\sigma(\overline f)$. 
Let $A\in C^g_n$ and define $\eta_{A,\sigma}\in {\mathfrak W}_{g-n}\times 
\widehat{\mathfrak W}_n$ by 
$\sigma^{-1}\pi_A=\pi_{\sigma^{-1}(A)}\eta_{A,\sigma}^{-1}$. Hence we have 
$\overline f'(A)=f'(\pi_A)=(\sigma f)(\pi_A)=f(\sigma^{-1}\pi_A)
= f(\pi_{\sigma^{-1}(A)}\eta_{A,\sigma}^{-1})=\eta_{A,\sigma}f(\pi_{\sigma^{-1}(A)})$. 
Thus 
$$
\begin{array}{lllll}
\Upsilon(\overline f') &=&\displaystyle 
\bigoplus_{A\in C^g_n}\pi_A(\overline f'(A))&=&\displaystyle 
\bigoplus_{A\in C^g_n}\pi_A\eta_{A,\sigma}\overline f(\sigma^{-1}A)=\\
&=&\displaystyle 
\bigoplus_{A\in C^g_n} \sigma\pi_{\sigma^{-1}(A)}\overline f(\sigma^{-1}A)
&=& 
\displaystyle 
\bigoplus_{B\in C^g_n} \sigma\pi_{B}\overline f(B)=\\
&=&
\displaystyle 
\sigma(\bigoplus_{B\in C^g_n}\pi_B(\overline f(B)))&=&
\displaystyle 
\sigma (\Upsilon(\overline f))\;, 
\end{array}
$$
which is what we needed to show.
Isometry of $\Upsilon$ is with respect to 
the natural inner product on 
 $Map(C_n^g,\, M^{g-n}_{\Z}\otimes L^n_{\Z})$  given by 
\beq\label{eq-indinner}
\lz\overline f, \overline h\rz\;\;=\;\;\sum_{A\in C^g_n} \lz \overline f (A), \overline g (A)\rz
\eeq
given the inner form on $M^{g-n}_{\Z}\otimes L^n_{\Z}$.  Also, as  $L^n_{\Z}$ is an
$\sl(2,\Z)$-module also $Map(C_n^g,\, M^{g-n}_{\Z}\otimes L^n_{\Z})$ is.
Both, equivariance and isometry, follows immediately from the form of 
the isomorphism in (\ref{eq-isomind}) and the fact that the $\pi_A$ are 
isometric equivariant maps. 
\ep
\smallskip

The relation of the isomorphisms in Lemma~\ref{lm-induce} to the ones in 
 Lemma~\ref{lm-upsilambda} is given by the restrictions of $\Upsilon^{-1}$
to the weight spaces 
\beq\label{eq-relupsis}
\pi_{N(\lambda)}^{-1}\;=\;\Upsilon_{\lambda^{\pi}}\circ\Upsilon_{\lambda}^{-1}
\;:\;{\cal W}(\lambda,g)\longrightarrow {\cal W}(\lambda^{\pi},g)\;. 
\;\;\subset \;\; {\cal W}(N_n,g)
\eeq
Here $\lambda^{\pi}$ denotes the the ``sign-content'' of 
a weight $\lambda\in\nabla_g$ defined by 
$$
\lambda^{\pi}:=\pi_{N(\lambda)}^{-1}\lambda=
\;\sum_{j=1}^{g-n(\lambda)}
\lambda_{\pi_{N(\lambda)}(j)}\epsilon_j
\;\;=\;\; \pm\epsilon_1+ 
\dots
\pm\epsilon_{g-n(\lambda)}\;. 
$$ 
\bigskip

\head{3. Lefschetz Decompositions and Specht Modules}\label{S3} 
\nopagebreak

As in \cite{Ker01} we consider the decomposition of the Frohman Nicas
TQFT according  to $SL(2,\R)$-representations:
$$
{\cal V}_{\Z}^{FN}\;=\;\;\bigoplus_{j\geq 1} V_{j}\otimes {\cal V}^{(j)}\;,   
$$
where $V_{j}$ is the $j$-dimensional irreducible representation
of $SL(2,\R)$. Note, that we the convention we use here for the 
superscript in ${\cal V}^{(j)}$ is shifted by one from the one used in 
\cite{Ker01}.  For weights we follow the notations of 
\cite{GooWal98}. The sublattices of the 
irreducible TQFT components can be defined as  the $SL(2,\R)$
lowest weight spaces
$$
{\cal V}^{(j)}_{\Z}(\Sigma)\;\;=\;\;\{v\in {\cal V}_{\Z}(\Sigma)\,:
\; Fv=0\;\;\;\mbox{and}\;\;\;Hv=-(j-1)v\}\;\;
=\;\;\ext {g-j+1}H_1(\Sigma,\Z)\cap ker(F)\;\;. 
$$
The representation of ${\rm Sp}(2g,\Z)$  on 
${\cal V}^{(j)}_{\Z}(\Sigma)$ is irreducible of fundamental 
heighest weight $\varpi_{g-j+1}=\epsilon_1+\ldots +\epsilon_{g-j+1}\,$
with heighets weight vector 
$$
w(\varpi_{g-j+1})=a_1\wedge\ldots\wedge a_{g-j+1}=\Upsilon_{\varpi_{g-j+1}}
(e_-^{\otimes j+1})\;. 
$$
The possible weights in this representation are given by 
$$
\nabla_{g}^{(j)}\;\;=\;\;\{\lambda\in\nabla_g\,:\;n(\lambda)\geq j-1\;\mbox{and}\;n(\lambda)\equiv j -1\mod 2\,\}\;.
$$
We obtain an analogous weight space decomposition
\beq\label{eq-weightcomp}
{\cal V}_{\Z}^{(j)}(\Sigma_g)\;\;=\;\;\bigoplus_{\lambda\in\nabla^{(j)}_g}
{\cal W}^{(j)}_{\Z}(\lambda,g)
\;\stackrel{\Upsilon\;\cong}{-\!\!\!-\!\!\!-\!\!\!-\!\!\!\longrightarrow}
\;\bigoplus_{\lambda\in\nabla^{(j)}_g}
L^{n(\lambda)}_{\Z}\cap {ker(F)}\cap ker(H+j-1)\;\;. 
\eeq

\begin{lemma}\label{lm-spechtiso}
The spaces ${\cal W}^{(j)}_{\Z}(\lambda,g)$ are as 
$S_{n(\lambda)}$-modules isomorphic to the standard irreducible
Specht modules ${\cal S}^{[a,b]}$ for the two-row Young-diagram 
$$
[a,b]\;=\;\Bigr [ \frac {n(\lambda)+j-1}2, \frac {n(\lambda)-j+1}2\Bigr]\;. 
$$
\end{lemma}

{\em Proof:} Although this appears to be standard we shall provide a proof
to fix conventions. 
We largely follow here the definitions and notations of \cite{Jam78}. 
First we note that $\,ker(H+j-1)$ is naturally isomorphic to
the permutation module $M^{[a,b]}$, where $n(\lambda)=a+b$ and
$j=a-b+1$. The isomorphism $M^{[a,b]}\stackrel{\cong}{\longrightarrow}
ker(H+j-1)\cap L^{n(\lambda)}_{\Z}\,$ maps a tabloid 
$\{t\}=\!\!\!\!\raise-8pt\hbox{
\bstid 
i_1 & i_2 & \ldots & i_a \cr
j_1& \ldots & j_b \cr
\estid}\,\,
$
to the basis vector  
$e_{\{t\}}\,=\,e_{\epsilon_1}\otimes e_{\epsilon_2}\otimes\ldots\otimes e_{\epsilon_N}$ with 
$\epsilon_k=+$ if $k\in\{j_1,\ldots,j_b\}$ and 
$\epsilon_k=-$ if $k\in\{i_1,\ldots,i_a\}$. It is obvious that the 
$e_{\{t\}}$ indeed span  ${\rm ker}(H+j-1)$.
Consider a tableau
$t=
\trA{i_1  \cr j_1 }\!
\trB{i_b & i_{b+1} \cr j_b }\!
\trC{i_a}\,
$
of shape $[a,b]$. Let $C_t$ be the column stabilizer group and 
\beq\label{eq-defkappa}
\kappa_t=\sum_{\pi\in C_t}sign(\pi)\pi=\prod_{k=1}^b(1-(i_k,j_k))
\eeq
the signed column sum. We set
$e_t\;=\;\kappa_t e_{\{t\}}\;$.
The Specht module ${\cal S}^{[a,b]}$
is the space generated by these $e_t$ and we want to show that
it coincider with $ker(F)$. The easy part is to show 
${\cal S}^{[a,b]}\subset ker(F)$. 
Using that $F$ commutes with $\kappa_t$ 
we compute $Fe_t=F\kappa_t e_{\{t\}}=\kappa_t Fe_{\{t\}}=
\sum_{k=1}^b\kappa_t e_{\{t_k\}}$, where $\{t_k\}$ is the tabloid of shape 
$[a+1,b-1]$, in which we have removed $j_k$ from the bottom row and added to 
the top row. As a result
$(1-(i_k,j_k))\{t_k\}=0$, hence $\kappa_t\{t_k\}=0$ so that $Fe_t=0$.

In order to prove $ker(F)\subset {\cal S}^{[a,b]}$ we proceed by induction.   $ker(F)$
on $M^{[a,b]}$ is given by $ker(F^2)$ on $M^{[a-1,b]}$. 
We have a map
$$
ker(F^2)\cap M^{[a-1,b]}\;\longrightarrow\; ker(F)\cap M^{[a,b]}\;\;:\;\;\;
x\;\mapsto\;e_-\otimes x\,-\,e_+\otimes Fx\;\;.
$$
It is easy to see that this map is an isomorphism. We use first that 
$Fx\in ker(F)\cap M^{[a,b-1]}$ and hence by induction $Fx=\sum_{t\in T} b_te_t$.
Here $T$ is the set of tableaux of the form 
$t=
\trA{i_1  \cr j_1 }\!
\trB{i_{b-1} & i_{b} \cr j_{b-1} }\!
\trC{i_a}\,
$
with numbers from $\{2,3,\ldots, a+b\}$. Now as $b-1<a$ we have that $i_a$ is not permuted
by $C_t$. Thus if we denote $E_j=\id^{\otimes j-1}\otimes E\otimes \id^{\otimes N-j}$ we have
that $E_{i_a}$ commutes with $\kappa_t$. Let
$
h_t=\kappa_tE_{i_a}e_{\{t\}}
$
Now $Fh_t=FE_{i_a}\kappa_t e_{\{t\}}=
[F,E_{i_a}]\kappa_t e_{\{t\}}+E_{i_a}F\kappa_t e_{\{t\}}
= -H_{i_a}\kappa_t e_{\{t\}}$ since $e_t\in ker(F)$ as shown above. Now,
$H_{i_a}$ commutes with $\kappa_t$ and $H_{i_a}e_{\{t\}}=-e_{\{t\}}$ by
construction. Hence
$
Fh_t=e_t
$
In other words ${\cal S}^{[a,b-1]}\subset im(F)$. Consider now $y=x-\sum_{t\in T}b_th_t$.
We thus have $y\in ker(F)\cap M^{[a-1,b]}$ so that by induction $y=\sum_{s\in S}c_s e_s$,
where $S$ denotes the tableau of the form
$
s=\!\!\!
\trA{k_1 \cr l_1}\!
\trB{k_b & k_{b+1}\cr l_b}\!
\trC{k_{a-1}}\,
$
with all numbers in $\{2,3,\ldots, a+b\}$. Inserting everything we find
$$
z\;=\;e_-\otimes x\,-\,e_+\otimes Fx\;\;=\;\; \sum_{s\in S}c_s e_-\otimes e_s\;\;+\;\;
 \sum_{t\in T} b_t(e_-\otimes h_t-e_+\otimes e_t)\;. 
$$
Now,  it is not hard to see that $e_{\hat s}= (e_-\otimes e_s)$ with
$
\hat s=\!\!\!
\trA{k_1 \cr l_1}\!
\trB{k_b & k_{b+1}\cr l_b}\!
\trC{k_{a-1} & 1}\,
$
given that $a-1\geq b$ so that 1 is not permuted by $C_{\hat s}$. Moreover,
$(e_-\otimes h_t-e_+\otimes e_t)=(1-(1,i_a))(e_-\otimes h_t)
(1-(1,i_a))\kappa_tE_{i_a}e_{\{t\}}=\kappa_{\hat t}e_{\{\hat t\}}=e_{\hat t}$, 
where
$
\hat t=\!\!
\trA{1 & i_1 \cr i_a & j_1}\!
\trB{i_{b-1} & i_b\cr j_{b-1}}\!
\trC{i_{a-1}}
$
Thus $z\;=\;\sum_{\hat s}c_se_{\hat s}\,+\,\sum_{\hat t}b_te_{\hat t}$
so that $z\in {\cal S}^{[a,b]}$. 
\ep
\smallskip 

Combining Lemma~\ref{lm-induce}  with Lemma~\ref{lm-spechtiso} and using 
notation ${\cal W}^{(j)}_{\Z}(A,g)$ and ${\cal W}^{(j)}_{\Z}(n,g)$ 
with $n\geq j-1$ and $ n \equiv j-1 \mod 2$ analogous to (\ref{eq-Wnotation})
so that  ${\cal V}_{\Z}^{(j)}(\Sigma_g)=\bigoplus_n
{\cal W}^{(j)}_{\Z}(n,g)$
 we find the following structure.

\medskip
 
\begin{coro}\label{cor-UpsSpecWei} 
For every $1\leq j-1\leq n\leq g$ with $n \equiv j-1 \mod 2$ 
there is an isomorphism 
of $\widehat{\mathfrak W}_g$-modules
$$
\Upsilon^{(j)}\;\;:\;\;\;
Ind^{\,\,\widehat{\mathfrak W}_g}_{\widehat{\mathfrak W}_{g-n}\times \widehat{\mathfrak W}_n}
\Bigl(M^{g-n}_{\Z}\otimes {\cal S}^{[\frac{n+j-1}2,\frac{n-j+1}2 ]}\Bigr)
\;\;\;\stackrel{\cong}{\longrightarrow}\;\;\;\;{\cal W}_{\Z}^{(j)}(n,g)\;. 
$$
\end{coro} 

\medskip

Let us also describe the $\sp (2g,\Z)$-generators on the vectors
$e_t$ spanning the Specht modules ${\cal W}^{(j)}(\lambda,g)$. 
To this end it is convenient to 
use tableaux in which entries are takes from the set $N(\lambda)$ rather
than $\{1,\ldots, n(\lambda)\}$, related to the standard ones by 
application of 
$\pi_{N(\lambda)}$. We denote the set of these tableaux 
by $T^{(j)}(\lambda)$.   

\bigskip

\begin{lemma}\label{lm-spechtTQFTgen}
The $\sp (2g,\Z)$-generators act on the  tableau vectors of Specht modules
in the ${\cal V}^{(j)}_{\Z}$-TQFT as follows. For the 
$e_{\alpha_i,\lambda}$ with $1\leq i\leq g-1$ we have: 
\begin{enumerate}
\item If $(\lambda_i,\lambda_{i+1})=(0,1)$ then 
$e_{\alpha_i,\lambda}e_t=e_s$ where $s\in T^{(j)}(\lambda+\alpha_i)$ is
obtained from $t\in T^{(j)}(\lambda)$ by replacing the label $i\in N(\lambda)$
in $t$ by the label ${i+1}\in N(\lambda+\alpha_i)$.
\item If $(\lambda_i,\lambda_{i+1})=(-1,0)$ then 
$e_{\alpha_i,\lambda}e_t=-e_s$ where $s\in T^{(j)}(\lambda+\alpha_i)$ is
obtained from $t\in T^{(j)}(\lambda)$ by replacing the label $i+1\in N(\lambda)$
in $t$ by the label $i\in N(\lambda+\alpha_i)$.
\item If $(\lambda_i,\lambda_{i+1})=(-1,1)$ we have 
$\{i,i+1\}= N(\lambda+\alpha_i)- N(\lambda)$ and 
$e_{\alpha_i,\lambda}e_t=e_s$ where $s\in T^{(j)}(\lambda+\alpha_i)$
 by adding a column 
\raise-7pt\hbox{$\bsyg i\cr i+1\cr \esyg$} to $t\in T^{(j)}(\lambda)$
\item If $(\lambda_i,\lambda_{i+1})=(0,0)$ so that 
$\{i,i+1\}= N(\lambda)- N(\lambda+\alpha_i)$ then 
$e_{\alpha_i,\lambda}e_t$, with  $t\in T^{(j)}(\lambda)$, is 
\begin{enumerate}
\item 0 if the labels $i$ and $i+1$ occur in columns of height 1. 
\item $2 e_s$ if $t$ is given by adding the column 
\raise-7pt\hbox{$\bsyg i\cr i+1\cr \esyg$} to $s$. 
\item $e_s$ if $i$ and $i+1$ occur in different columns of height 2. 

Here $s$
is obtained from $t$ by replacing the double column 
\raise-7pt\hbox{$\bsyg k & l \cr i & i+1\cr \esyg$} by
\raise-7pt\hbox{$\bsyg l \cr k \cr \esyg$}
\item $e_s$ if $i$ is in column of height 2 and and $i+1$ in 
columns of height 1, where $s$
is obtained from $t$ by deleting the column 
\raise-7pt\hbox{$\bsyg i \cr k \cr \esyg$} and 
replacing \raise-3pt\hbox{$\bsyg i +1  \cr \esyg$} by 
\raise-3pt\hbox{$\bsyg k \cr \esyg$}
\item $e_s$ if $i+1$ is in column of height 2 and and $i$ in 
columns of height 1, where $s$
is obtained from $t$ by deleting the column 
\raise-7pt\hbox{$\bsyg k \cr i+1 \cr \esyg$} and 
replacing 
\raise-3pt\hbox{$\bsyg i  \cr \esyg$} by 
\raise-3pt\hbox{$\bsyg k \cr \esyg$}
\end{enumerate}
\end{enumerate}
All other cases of positions of $i$ and $i+1$ follow from the 
symmetry properties of the vectors $e_t$ under permutations of columns or 
within columns. Since $N(\lambda+\alpha_g)=N(\lambda)$ if 
$\lambda+\alpha_g, \lambda\in \nabla_g$ we have that $e_{\alpha_g, \lambda}$
acts as identity also on the vectors $e_t$. Similarly, $H_{\lambda}$ acts as
identity. 
\end{lemma}

\bigskip

\head{4. $\F_p=\Z/p\Z$-Reductions and the Sequences  ${\cal C}_{p,k}$}\label{S4}
\nopagebreak

For the remainder of this article let $p$ be an odd prime number.
Since the TQFT's  ${\cal V}^{(j)}_{\Z}$ are defined over 
free $\Z$-modules (lattices) we naturally obtain TQFT's ${\cal V}^{(j)}_p$ 
over the number field $\F_p$ by setting
\vspace*{-6mm}

$$
{\cal V}^{(j)}_p(\Sigma_g)\;=\;
\raise5pt\hbox{${\cal V}^{(j)}_{\Z}(\Sigma_g)$}\Big / 
\raise-5pt\hbox{$p {\cal V}^{(j)}_{\Z}(\Sigma_g)$}
$$
Now,  each ${\cal V}^{(j)}_{\Z}(\Sigma_g)$
inherits a non-degenerate inner product as a  sublattice of
$\ext *H_1(\Sigma_g,\Z)$. 
This, however, will in general degenerate if we consider the 
the $p$-reduction 
$\lz\,,\,\rz_p:\Bigr({\cal V}^{(j)}_p(\Sigma_g)\Bigl)^{\otimes 2}\to \F_p$.  
We denote the corresponding null space as follows. 
$$
\nil {\cal V}^{(j)}_p(\Sigma_g)\;\;=\;\;\bigl\{v\in {\cal V}^{(j)}_p(\Sigma_g)\,:\;\lz v,w\rz_p=0\;\;
\forall\,w\in\,{\cal V}^{(j)}_p(\Sigma_g) \bigr\}
$$
The elements are represented by vectors $v\in {\cal V}^{(j)}_{\Z}(\Sigma_g)$ for which $\lz v,w\rz\in p\Z$
for all $w\in {\cal V}^{(j)}_{\Z}(\Sigma_g)$ although $v\not \in p{\cal V}^{(j)}_{\Z}(\Sigma_g)$. 
\begin{lemma} There are well defined TQFT's $\nil {\cal V}^{(j)}_p$ and $\dov {\cal V}^{(j)}_p$
which assign to a surface $\Sigma_g$ the $\F_p$-vectors spaces  
\vspace*{-6mm}

$$
\qquad\qquad \nil {\cal V}^{(j)}_p(\Sigma_g)\qquad\qquad\mbox{and}\qquad\qquad 
\dov {\cal V}^{(j)}_p(\Sigma_g)\;=\;
\raise5pt\hbox{${\cal V}^{(j)}_p(\Sigma_g)$}\Big / 
\raise-5pt\hbox{$\nil {\cal V}^{(j)}_p(\Sigma_g)$}
$$
\end{lemma}

{\em Proof:} We note that since for a cobordism $M$ the map
${\cal V}_{\Z}(M)$ commutes with $E$ we have that
${\cal V}_{\Z}(M)^*$ commutes with $F=E^*$ and hence also maps the ${\cal V}_{\Z}^{(j)}(\Sigma_g)$ to 
themselves. Thus if $v_i\in{\cal V}^{(j)}_{\Z}(\Sigma_{i})$ for $i=1,2$, $M$ is a cobordism from
$\Sigma_{g_1}$ to $\Sigma_{g_2}$, and $v_1$ represents a vector in 
$\nil {\cal V}^{(j)}_p(\Sigma_{1})$, then $\lz v_2,{\cal V}^{(j)}_{\Z}(M)v_1\rz=
\lz {\cal V}^{(j)}_{\Z}(M)^*v_2,v_1\rz\in p\Z$ as 
${\cal V}^{(j)}_{\Z}(M)^*v_2\in {\cal V}^{(j)}_{\Z}(\Sigma_{1})$.
\ep
\smallskip

We extend the previous notations to the weight spaces $\nil{\cal W}^{(j)}_p(\lambda,g)$,  
$\dov {\cal W}^{(j)}_p(\lambda,g)$, $\nil{\cal W}^{(j)}_p(n,g)$, etc.. Since the weight spaces
are all orthogonal to each other these subspaces can be defined also as the null spaces from 
the respective  restriction of the inner forms. Now, also the Specht modules ${\cal S}^{\tau}$
for a diagram $\tau=[a,b]$ inherit an inner form from the permutation module 
$M^{\tau}$, which is via the isometry $\Upsilon$ compatible with the one on the weight spaces. 
As in the standard literature, e.g.,  \cite{Jam78}, we set 
$$
{\cal D}^{\tau}_p=\raise4pt\hbox{${\cal S}^{\tau}_p$}\Big / 
\raise-4pt\hbox{$\nil {\cal S}^{\tau}_p$}\qquad\quad\mbox{where}\quad 
\nil {\cal S}^{\tau}_p={\cal S}^{\tau}_p\cap {{\cal S}^{\tau}_p}^{\perp}
$$
and ${\cal S}^{\tau}_p$ is the $p$-reduction of ${\cal S}^{\tau}$. They are related to
irreducible  TQFT's as follows. 

\begin{lemma} Let $p\geq 3$ be a prime. 
The TQFT's $\dov {\cal V}^{(j)}_p$ are {\em irreducible} over $\F_p$ and 
the weight spaces are identified as 
$\widehat{\mathfrak W}_g$-modules by equivariant isomorphisms: 
\beq\label{eq-weightiso-p}
\Upsilon^{(j)}_p\;\;:\;\;\;
Ind^{\,\,\widehat{\mathfrak W}_g}_{\widehat{\mathfrak W}_{g-n}\times \widehat{\mathfrak W}_n}
\Bigl(M^{g-n}_p\otimes {\cal D}^{[\frac{n+j-1}2,\frac{n-j+1}2 ]}_p\Bigr)
\;\;\;\stackrel{\cong}{\longrightarrow}\;\;\;\;\dov {\cal W}_p^{(j)}(n,g)\;. 
\eeq
\end{lemma}

{\em Proof:} The isomorphism in (\ref{eq-weightiso-p}) follows from the definitions and properties
of $\Upsilon^{(j)}$. We first show that the spaces $\dov {\cal W}_p^{(j)}(n,g)$ are irreducible with 
respect to the semidirect product ${\mathfrak X}_g$ of the Cartan algebra 
$\Z[{\mathfrak h}_g]\subset \sl(2g,\Z)$ and the algebra of the Weyl group 
$\Z[\widehat{\mathfrak W}_g]\subset \Z[{\rm Sp}(2g,\Z)]$.
For any vector $v=\sum_{\lambda} v_{\lambda}\in \dov {\cal W}_p^{(j)}(n,g)$ the action of 
${\mathfrak h}_g$ shows that each weight component $v_{\lambda}\in \dov {\cal W}_p^{(j)}(\lambda,g)$
has to lie in ${\mathfrak X}_gv$. Each $\dov {\cal W}_p^{(j)}(n,g)$ is a module 
of a 
symmetric group $S_n\subset \widehat{\mathfrak W}_g\subset{\mathfrak X}_g$ which is 
equivalent to ${\cal D}^{[\frac{n+j-1}2,\frac{n-j+1}2 ]}_p$. It  now follows from Theorem~4.9 in \cite{Jam78}
that these representations are irreducible. In fact, as $p\geq 3$ any two-row diagram is $p$-regular
so that these representations are never zero, see Theorem~11.1 in \cite{Jam78}. In particular, if 
$v_{\lambda}\neq 0$ then $\Z[S_n]v_{\lambda}$ is the entire module. Hence, for $v\neq 0$ we
must have $\dov{\cal W}_p^{(j)}(\lambda,g)\subset {\mathfrak X}_gv$ for at least one $\lambda\in\nabla_g$ with
$n(\lambda)=n$. Since $\widehat{\mathfrak W}_g$ acts transitively on all of such weights and 
provides isomorphisms between the weight spaces we 
thus have $\dov {\cal W}_p^{(j)}(n,g)= {\mathfrak X}_gv$, which implies irreducibility since $v$ was 
arbitrary. A submodule of $\dov {\cal V}_p^{(j)}(\Sigma_g)$ must therefore be a direct sum of the
$\dov {\cal W}_p^{(j)}(n,g)$. Each of these contains a special vector 
$w_n^g=\Upsilon_{\varpi_{g-j+1}}(e_{t(n,j)})$ with  
 $t(n,j)=
\trA{1 & 3 \cr 2 & 4 }\!
\trB{m-1 & m+1 \cr m  }\!
\trC{n}\,
$ where $m=n-j+1$.  Using 
$\kappa_{t(n,j)}^2=2^m\kappa_{t(n,j)}$ for the antisymmetrizer (\ref{eq-defkappa}),
we find $\lz w^g_n, w^g_n\rz_p=\lz e_{t(n,j)}, e_{t(n,j)}\rz_p=
2^m\lz e_{\{t(n,j)\}},\kappa_t e_{\{t(n,j)\}}\rz_p=2^m\not\equiv 0\mod p$ for $p\geq 3$
so that all of these vectors are non-zero in  $\dov {\cal W}_p^{(j)}(n,g)$. 
It can be computed from the rules {\em 3.} and {\em 4.(b)} 
in Lemma~\ref{lm-spechtTQFTgen} that 
$e_{\alpha_{g-n-1}}{\sf S}_{g-n-1}w^g_n=w^g_{n-2}$ and $e_{\alpha_{g-n+1}}w^g_n=2w^g_{n+2}$,
where 
${\sf S}_l\in{\mathfrak W}_g$ maps $a_j$ to $-b_j$ so that with $2\not\equiv 0 \mod p$ 
we have non-trivial maps between all of the irreducible ${\mathfrak X}_g$-components.
Consequently, the $\dov {\cal V}_p^{(j)}(\Sigma_g)$ are irreducible as 
$\sp(2g,\Z)$-representations and any sub-TQFT must assign either this space or 0
to a surface $\Sigma_g$. As in \cite{Ker01} we easily check that the handle attachment
maps are non-trivial between these spaces so that $\dov {\cal V}_p^{(j)}$ does in fact
contain no proper sub-TQFT. 
\ep 
\smallskip

We next construct a sequence of maps between the $p$-reductions of the Specht modules,
using the $\sl(2,\Z)$-actions. As before we fix $n\in\N$  and 
denote the Specht module  
$$
{\cal S}^{\{c\}}_{\Z}\;\;=\;\;L_{\Z}^n\cap\ker(F)\cap\ker(H+c-1)\;\;\cong\;\;{\cal S}^{\tau}\;, 
$$
where $c=a-b+1$ with tableau $\tau=[a,b]=[\frac {n+c-1}2, \frac {n-c+1}2]$.

\begin{lemma}\label{lm-Emaps} Let  $p\geq 3$, $c\not\equiv 0\mod p$, and $c_0\in\{1,\ldots,p-1\}$
such that $c\equiv c_0\mod p$. Then 
$$
E^{c_0}(S^{\{c\}}_{\Z})\;\subset\;S^{\{c-2c_0\}}_{\Z}\,+\,pL^n_{\Z}\;\;.
$$
Moreover, 
$$
E^{c_0}(S^{\{c\}}_{\Z})\,\not\subset\,pL^n_{\Z}\qquad\quad
\mbox{and}\qquad\quad E^{p}(L^n_{\Z})\subset pL^n_{\Z}\;.
$$
\end{lemma}

{\em Proof:} Now, with notation as in the proof of Lemma~\ref{eq-weightiso-p}, 
 $e_{t(n,c-1)}\in S^{\{c\}}_{\Z}$ is a cyclic vector so it suffices to
show that $E^{c_0}(e_{t(n,c-1)})\in (S^{\{c-2c_0\}}_{\Z}\,+\,pL^n_{\Z})-pL^n_{\Z}$.
Furthermore, $e_{t(n,c-1)}=e_{t(m,0)}\otimes e_-^{\otimes c-1}$ and $Ee_{t(m,0)}=0$
so that we really need to show that 
$E^{c_0}(e_-^{\otimes c-1})\in (ker(F)\,+\,pL^{c-1}_{\Z})-pL^{c-1}_{\Z}$,
where $c=c_0+kp$. We do this by induction in $c_0$.

For $c_0=1$ we have $c-1=kp$ and, using $\sl_2$-relations, $FE(e_-^{\otimes kp})= 
-He_-^{\otimes kp}+EFe_-^{\otimes kp}=kp e_-^{\otimes kp}$. Hence, if we set 
$w=k e_-^{\otimes kp-1}\otimes e_+$ we find $Fw=e_-^{\otimes kp}$ so that
$E(e_-^{\otimes kp})-kpw\in ker(F)$, meaning $E(e_-^{\otimes kp})\in  ker(F)+pL^{c-1}_{\Z}\,$.
Next assume the assertion is true for $c_0$ so that 
$E^{c_0}(e_-^{\otimes c-1})=y\,+\,pz\,$ with $Fy=0$ and $-Hy=(c-2c_0-1)y$. 
. We have by binomial formula and $Ee_-=e_+$
and $E^2e_-=0$ that 
$E^{c_0+1}(e_-^{\otimes c})=E^{c_0+1}(e_-^{\otimes c-1}\otimes e_-)=
(E^{c_0+1}e_-^{\otimes c-1})\otimes e_-\,+\,(c_0+1)
E^{c_0}e_-^{\otimes c-1})\otimes e_+=(Ey)\otimes e_-\,+\,(c_0+1)y\otimes e_+\,+\,pz'\,$. 
Thus we need to show $t=(Ey)\otimes e_-\,+\,(c_0+1)y\otimes e_+\in ker(F)\cap pL^{c}_{\Z}$.
We compute $Ft=(FEy)\otimes e_-\,+\,(c_0+1)y\otimes e_-=
(-Hy)\otimes e_-\,+\,(c_0+1)y\otimes e_-=(c-c_0)y\otimes e_-=kp y\otimes e_-$. 
Also $F(y\otimes e_+)=y\otimes e_-$ so that 
$t-kp y\otimes e+-\in ker(F)$ and hence $t\in  ker(F)\cap pL^{c}_{\Z}$.

Finally, it is not hard to see that 
$\lz e_+^{\otimes c_0}\otimes e_-^{\otimes c-c_0-1}, E^{c_0}(e_-^{\otimes c-1})\rz = 
c_0 !\not\equiv 0\mod p$ if $c_0<p$ so that $E^{c_0}(e_-^{\otimes c-1})\not\in pL^{c-1}_{\Z}$
and hence $E^{c_0}e_{t(n,c-1)}\not\in pL^{n}_{\Z}$. Also, the binomial formula yields
$E^p=\sum_j{p\choose j}E^{p-j}\otimes E^j$ so that we can conclude $E^p=0$ by a similar 
induction argument.
\ep
\smallskip

In particular Lemma~\ref{lm-Emaps} implies that we have well defined,
non-zero maps 
$E^{c_0}\,:\,S^{\{c\}}_p\,\longrightarrow\, S^{\{c-2c_0\}}_p$ on the respective 
$p$-reductions $S^{\{c\}}_p=S^{\{c\}}_{\Z}/pS^{\{c\}}_{\Z}=
S^{\{c\}}_{\Z}/pL^n_{\Z}$ using that $pS^{\{c\}}_{\Z}=S^{\{c\}}_{\Z}\cap pL^n_{\Z}$. 

\begin{coro}\label{cor-Cseq}
 For $p$ and $n$ as above and $k=1,\ldots,p-1$ with $k\equiv n+1\mod 2$
there is a sequence  ${\cal C}_{p,k}$ 
of Specht modules  over $\F_p$ as follows: 
\beq\label{eq-Cseq}
0\,\to\,{\cal S}^{\{n+1-2l\}}_p\,\to\;\ldots\;\to\,
{\cal S}^{\{2p+k\}}_p\,\stackrel{E^{k}}{-\!\!\!-\!\!\!\longrightarrow}\,
{\cal S}^{\{2p-k\}}_p\,\stackrel{E^{p-k}}{-\!\!\!-\!\!\!-\!\!\!\longrightarrow}\,
{\cal S}^{\{k\}}_p\,\to\,{\cal D}^{\{k\}}_p\,\to\,0\;. 
\eeq
All maps (except the first and last one) 
are non-zero, and any two consecutive maps compose to zero.
\end{coro}

More precisely, we have that the $i$-th component of this sequence is 
 $\,{\cal C}^{(i)}_{p,k}={\cal S}^{\{ip+k_i\}}_p$, where $k_i=k$ if $i$ is even and 
$k_i=p-k$ if $i$ is odd. The maps are 
$E^{k_i}:{\cal C}^{(i)}_{p,k}\to {\cal C}^{(i-1)}_{p,k}$ so that two consecutive 
maps compose as  $E^{k_i}E^{k_{i+1}}=E^p=0$. 

We have that 
${\cal C}^{(0)}_{p,k}={\cal S}^{\{k\}}_p\to {\cal C}^{(-1)}_{p,k}={\cal D}^{\{k\}}_p$ 
is the (non-zero) quotient map.
Now, it is clear that $im(E)\subset ker(F)^{\perp}\subset \nil {\cal S}^{\{k\}}_p$
given that $E^*=F$. Hence also the composite 
${\cal C}^{(1)}_{p,k}\to{\cal C}^{(0)}_{p,k}\to {\cal C}^{(-1)}_{p,k}$ is zero. 
In order to characterize the last index write $\frac{n+1+k} 2= ph+q$ with $h\in\Z$ and 
$q=0,\ldots\, p-1$. We have 
\beq\label{eq-lastindex}
l=\left\{\begin{array}{ll} q-k\;&\;\mbox{if}\;q\geq k\\
                           q \;&\;\mbox{if}\;q<k\\
          \end{array}\right.\;. 
\eeq
 
The maps in Corollary~\ref{cor-Cseq} thus extend to 
 a sequence of the $p$-reductions of
 the induced representations from Corollary~\ref{cor-UpsSpecWei}
as well as  the weight spaces ${\cal W}_p^{(c)}(n,g)$. 
The respective maps on the vector spaces  ${\cal V}_p^{(c)}(\Sigma_g)$  
are, by equivariance of the $\Upsilon^{(j)}$, given by the restriction and
$p$-reductions of the maps $E^{c_0}$ on $\ext *H_1(\Sigma_g)$. In particular,  these
maps commute with the TQFT-images of the cobordisms by equivariance of ${\cal V}_{\Z}$.

\begin{coro}\label{cor-seqTQFT}
 For $p$, $n$  and $k$ as above 
there is a sequence of natural transformations of TQFT's 
\beq\label{eq-seqTQFT}
\ldots\,\to\,{\cal V}^{((i+1)p+k_{i+1})}_p\,\to\,
{\cal V}^{(ip+k_i)}_p\,\to\,\ldots\,\to\,
{\cal V}^{(2p-k)}_p\,\to\,
{\cal V}^{(k)}_p\,\to\,\dov {\cal V}^{(k)}_p\,\to\,0
\eeq
with $k_i=k, p-k$ as above and any two consecutive transformations compose to zero.
\end{coro}

\head{5. Exactness of   ${\cal C}_{p,k}$}\label{S5}
\nopagebreak

Exactness of the sequence  ${\cal C}_{p,k}$ defined in 
Corollary~\ref{cor-Cseq} follows from the modular structure 
of the involved Specht modules.  
The irreducible factors of ${\cal S}^{\tau}_p$ for two row diagrams $\tau$ are
determined in Theorem~24.15 of \cite{Jam78}. For our proof we will, however, have to
make use also of the precise submodule structure, which turns out to be rather
rigid. We use the recent result of Kleshchev and Sheth in \cite{KS99} that describes
this structure precisely. In order to state it we need to introduce some more 
conventions and definitions. 

For a Young diagram $\tau=[a,b]$  we set $c=a-b+1$, and consider the 
$p$-adic expansion 
\beq\label{eq-cpadic}
\qquad \qquad c\;=\;c_0\,+\,c_1p\,+\,c_2p^2\,+\ldots\,+\,c_rp^r\qquad\quad \mbox{with}
\quad c_j\in\{0,1,\ldots,p-1\}\;. 
\eeq
As in \cite{KS99} we denote by $\hat A_{\tau}$ the family of sets of integers of the 
form
\beq\label{eq-defIA}
I\,=\,[i_1,i_2)\cup [i_3,i_4)\cup \ldots\cup [i_{2u-1},i_{ut})
\eeq
such that
\vspace*{-8mm}

\beq 
i_1<i_2<\ldots < i_{2u}\qquad,\qquad c_{i_{2j-1}}\neq 0 \qquad\mbox{and}
\qquad c_{i_{2j}}\neq p-1. 
\eeq
For such a set $I\in A_{\tau}$ we define as in \cite{KS99}the number 
\beq\label{eq-defdelta} 
\delta^{\tau}_I\;=\;\sum_{i\in I} (p-1-c_i)p^i \;+\;\sum_{j=1}^up^{i_{2j-1}}
\;\;=\;\;\sum_{j=1}^u\delta^{\tau}_{[i_{2j-1},i_{2j})}\;,
\eeq
 where \ 
$\displaystyle\;\delta_{[u,w)}^{\tau}=\sum_{i=u}^{w-1}(p-1-c_i)p^i\,+\,p^u\,$. 
Also as in \cite{KS99} we introduce the smaller set $A_{\tau}\subset \hat A_{\tau}$ given by 
\beq\label{eq-defAb} 
A_{\tau}=\{I\in \hat A_{\tau}\,:\;\delta^{\tau}_I\leq b\}\;, 
\eeq
as well as for any $I\in A_{\tau}$
the function $\nu_I$ on Young diagrams $\tau=[a,b]$ defined as 
\beq\label{eq-defnu}
{\nu}_I(\tau)=[a+\delta_I^{\tau}, b-\delta_I^{\tau}]\;. 
\eeq
As before we denote by ${\cal D}^{\mu}_p$ the irreducible quotient of the Specht module 
${\cal S}^{\mu}_p$ over $\F_p$ for any two row diagram $\mu$. Also denote by 
${\cal F}(M)$ the set of irreducible factors that occur in a 
composition series of a representation $M$. Further, let  
${\cal M}_{\mu}^{\lambda}$ be the smallest sub module of $S^{\lambda}_p$ such that
${\cal D}^{\mu}_p\in {\cal F}({\cal M}_{\mu}^{\lambda})$. The description of 
the submodule structure
in \cite{KS99} uses the partial order on ${\cal F}({\cal S}^{\tau}_p)$ defined as 
${\cal D}_p^{\mu_1}\leq_{\tau}{\cal D}_p^{\mu_2}$ if and only if 
${\cal M}_{{\mu}_1}^{\tau}\subseteq 
{\cal M}_{{\mu}_2}^{\tau}$, i.e., if and only if ${\cal D}_p^{\mu_1}\subseteq 
{\cal M}_{{\mu}_2}^{\tau}$. 
 
\noindent
\begin{theorem}[Corollary 3.4 of  \cite{KS99}] \label{thm-KS99}\ 
\begin{enumerate}
\item All multiplicities of ${\cal D}^{\mu}$ in $S^{\tau}$ are zero or  one. 
\item ${\cal F}({\cal S}_p^{\tau})=\{{\cal D}_p^{\nu_I(\tau)}\,:\;I\in A_{\tau}\}$
\item ${\cal D}_p^{\nu_J(\tau)}\geq_{\tau} {\cal D}_p^{\nu_I(\tau)}$
if and only if $J\subseteq I$. \ep
\end{enumerate}
\end{theorem}

\noindent 
We define now a submodules ${\cal S}_p^{\tau}$ by choosing a special subset of of $A_{\tau}$
and $\hat A_{\tau}$. It is defined as 
\beq\label{eq-defAsets}
\hat A^0_{\tau}=\{I\in \hat A_{\tau}\,:\; 0\in I\}
\qquad\quad
\mbox{and}\qquad\quad
\hat A^+_{\tau} = \hat A-\hat A^0_{\tau}=\{I\in \hat A_{\tau}\,:\; 0\not\in I\}\;
\eeq
as well as\ \ \qquad
$
A^0_{\tau}=\hat A^0_{\tau}\cap  A_{\tau}
\qquad\quad
\mbox{and}\qquad\quad
A^+_{\tau}=\hat A^+_{\tau}\cap  A_{\tau}
$.  
Let us also introduce the number $k_{\tau}= min\{j\geq 1:\;c_j\neq 0\}$ with $k_{\tau}=\infty$
if $c=c_0$. Hence
\beq\label{eq-Atauplus}
c\;=\;c_0\,+\,c_{k_{\tau}}p^{k_{\tau}}\,+\,\ldots\,+\,c_rp^r
\qquad\quad\mbox{and}
\qquad 
[0,k_{\tau})\cap I=\emptyset \qquad\mbox{for any}\quad I\in \hat A_{\tau}^+\;. 
\eeq
The latter follows since if $0\not\in I$ we must have $i_1>0$. But with $c_{i_1}\neq 0$
we find $i_1\geq k_{\tau}$. 
For the following we assume $c\not\equiv 0\mod p$, i.e., $c_0\neq 0$. Given this
we introduce for a diagram $\tau=[a,b]$ with $a-b=c-1\geq 2c_0$ 
(so that $k_{\tau}<\infty$) the notation
\beq\label{eq-tauprime}
\tau'\;=\;[a-c_0,b+c_0]\qquad\quad\mbox{so that}\qquad \tau=\nu_{[0,k_{\tau})}(\tau')
\eeq
For $\tau'$ we thus have $c'=c-2c_0$ and for 
the $p$-adic expansion of $c'=\sum_jc'_jp^j$ we obtain
\beq\label{eq-ccoeff}
\begin{array}{rcll}
c_0' &=& p-c_0 & \\ 
c_i' &=& p-1 & \quad\mbox{for}\; 1\leq i <k_{\tau}\\
\end{array}
\qquad\quad
\begin{array}{rcll}
c'_{k_{\tau}} &=& c_{k_{\tau}}-1&\\
c'_i\, &=& c_i &\;\mbox{for}\; i>{k_{\tau}}\\
\end{array}
\eeq
From these equations it is clear that $c_0'\neq 0$ and $c'_{k_{\tau}}\neq p-1$ so 
that $[0,k_{\tau})$ is an admissible interval for $\hat A_{\tau'}$. In fact, it is
the unique minimal interval of the special subset
\beq\label{eq-AprimeInt}
\hat A_{\tau'}^0\;=\;\{I\in \hat A_{\tau'}:\;[0,k_{\tau})\subseteq I\}\;
\qquad\mbox{and}\qquad [0,k_{\tau})\in A_{\tau'}^0\;. 
\eeq
This is obvious since we must have $i_1=0$ for $0\in I$ and then the next possible $i_2$ 
is $k_{\tau}$. It is also easy to see that $\delta_{[0,k_{\tau}]}^{\tau'}=c_0$,
confirming the relation in (\ref{eq-tauprime}) between $\tau$ and $\tau'$ via
(\ref{eq-defnu}). We have the following simple but crucial observation.

\begin{lemma}\label{lm-hatbij} 
With $\tau$ and $\tau'$ as above we have a well defined bijection with inverse
$$
\phi:\;\hat A^0_{\tau'}\,\longrightarrow\, \hat A^+_{\tau}\;\;
:\;I\;\mapsto\;I-[0,k_{\tau})\quad\quad\mbox{and}\quad\quad
\phi^{-1}:\;\hat A^+_{\tau}\,\longrightarrow\, \hat A^0_{\tau'}\;\;
:\;J\;\mapsto\;J\cup [0,k_{\tau})
$$
\end{lemma} 

{\em Proof:}  We first show that $\phi$ is well defined. From (\ref{eq-AprimeInt})
we know that every $I\in  A_{\tau'}^0$ contains the special interval and is thus of the
form 
$$
I=[0,i_2)\cup [i_3,i_4)\cup\ldots \cup [i_{2u-1}, i_{2u})\qquad\quad\mbox{with}
\quad i_2\geq k_{\tau}\;. 
$$
Clearly,  $0\not\in\phi(I)$. We further distinguish the following two cases:

\paragraph{ Case $i_2={k_{\tau}}$ :} In this situation $\phi(I)= [i_3,i_4)\cup\ldots\;$.
Now $\phi(I)\in \hat A^+_{\tau}$ since by (\ref{eq-ccoeff}) 
the coefficients of $c$ and $c'$
and hence the conditions on the $i_j$ for $j\geq 3$ are 
the same.

\paragraph{ Case $i_2>{k_{\tau}}$ :} In this situation $\phi(I)= 
[{k_{\tau}},i_2)\cup R$ with 
$R=\,[i_3,i_4)\cup\ldots\;$. 
We have $c_{k_{\tau}}\neq 0$ 
by assumption, and $c_{i_2}=c'_{i_2}\neq q-1$ as $i_2>{k_{\tau}}$.
Also all intervals in $R$ are admissible. Hence $\phi(I)\in \hat A^+_{\tau}$
again. 

In a similar fashion we show that $\phi^{-1}$ is well defined using the fact that
any $L\in \hat A^+_{\tau}$ is by (\ref{eq-Atauplus}) of the form 
$L=[j_1,j_2)\cup[j_3,j_4)\cup\ldots\,$ with $j_1\geq k_{\tau}$. For $j_1> k_{\tau}$
we have $\phi^{-1}(L)=[0,k_{\tau})\cup [j_1,j_2)\cup[j_3,j_4)\cup\ldots\,$. This is
$\hat A^0_{\tau'}$ as $c'_0\neq 0$, $c'_{k_{\tau}}\neq p-1$,  and $c'_j=c_j$ for
$j\geq j_1$. In case that $j_1= k_{\tau}$ we have 
$\phi^{-1}(L)=[0,j_2)\cup[j_3,j_4)\cup\ldots\,$, which is again in 
$\hat A^0_{\tau'}$ since $c'_j=c_j$ for
$j\geq j_2>j_1$.
It is obvious that $\phi$ and $\phi^{-1}$ are inverses of each other given
(\ref{eq-Atauplus}) and (\ref{eq-AprimeInt}). 
\ep 

Next we consider how the differential from (\ref{eq-defdelta}) changes under this bijection.
\begin{lemma}\label{lm-delta} We have 
$$
\delta^{\tau'}_I \;=\;\delta^{\tau}_{\phi(I)}\,+\,c_0
$$
\end{lemma}

\noindent
{\em Proof:} This is a computation done according to the same cases as in
the proof of Lemma~\ref{lm-hatbij}:
\paragraph{ Case $i_2=k_{\tau}$ :} Then 
$\delta_{I}^{\tau'}=\delta_{[0,{k_{\tau}})}^{\tau'}+\sum_{j\geq 2}
\delta^{\tau'}_{[i_{2j-1},i_{2j})}$. Now, since $c_i'=p-1$ for 
$i=1,\ldots,{k_{\tau}}-1$ we have $\delta_{[0,{k_{\tau}})}^{\tau'}=(p-1-c'_0)+p^0=c_0$. 
We also have $\delta^{\tau'}_{[i_{2j-1},i_{2j})}=\delta^{\tau}_{[i_{2j-1},i_{2j})}$
since the $c_i$ are all the same with $i_j>{k_{\tau}}$. But 
$\delta^{\tau}_{I-[0,{k_{\tau}})}= 
\sum_{j\geq 2} \delta^{\tau}_{[i_{2j-1},i_{2j})}$ which implies the assertion. 
\paragraph{ Case $i_2>k_{\tau}$ :} As before we write
  $I=[0,i_2)\cup[i_3,i_4)\cup\ldots\,=[0,i_2)\cup R$ and $\phi(I)= [k_{\tau},i_2)\cup R$
\begin{eqnarray}
\delta_{I}^{\tau'}&=& \sum_{i=0}^{i_2-1}(p-1-c_i')p^i + 1\;+\;\delta_R^{\tau'}\\
&=& (p-1-c'_0)+1 \,+\,(p-1-c_{k_{\tau}}')p^{k_{\tau}}\,+
\,\sum_{i>{k_{\tau}}}^{i_2-1}(p-1-c_i')p^i\;+\;\delta_R^{\tau}\\
&=& c_0\,+\,(p-1-c_{k_{\tau}})p^{k_{\tau}}\,+\,p^{k_{\tau}}\,+
\,\sum_{i>{k_{\tau}}}^{i_2}(p-1-c_i)p^i\;+\;\delta_R^{\tau}\\
&=& c_0\,+\, \sum_{i={k_{\tau}}}^{i_2-1}(p-1-c_i)p^i\,+\,p^{k_{\tau}}
 \;+\;\delta_R^{\tau}\\
&=& c_0\,+\, \delta_{[{k_{\tau}},i_2)}^{\tau}\;+\;\delta_R^{\tau}
\;\;\;\;\;\;=\;\;\;\;\;\;   c_0\,+\,\;\delta_{[{k_{\tau}},i_2)\cup R}^{\tau}
\;\;\;\;\;\;=\;\;\;\;\;\;
c_0\,+\,\;\delta_{I-[0,{k_{\tau}})}^{\lambda}
\end{eqnarray}
which proves the assertion.
\ep
\smallskip

\begin{coro}\label{cor-bij} The map $\phi$ restricts to a bijection
$
\phi:\; A^0_{\tau'}\,\longrightarrow\,  A^+_{\tau}\;\;\;
$
with the properties that it is monotonous with respect to inclusions 
and 
\begin{equation}\label{eq-nuform}
{\nu}_I(\tau')\;=\;\nu_{\phi(I)}(\tau)\;.
\end{equation}
\end{coro}

{\em Proof:} From (\ref{eq-tauprime}) we have $b'=b+c_0$. Moreover, 
Lemma~\ref{lm-hatbij} and Lemma~\ref{lm-delta} imply that 
$I\in A^0_{\tau'}$ iff $\delta^{\tau'}_I\leq b'$
iff $\delta^{\tau}_{\phi(I)}+c_0\leq b+c_0$ iff 
$\delta^{\tau}_{\phi(I)}\leq b$ iff $\phi(I)\in A^+_{\tau}$.
The fact that $\phi(A)\subset \phi(B)$ iff $A\subset B$ is obvious. 
Relation (\ref{eq-nuform}) is immediate from Lemma~\ref{lm-delta},
(\ref{eq-defnu}) and (\ref{eq-tauprime}).
\ep

Let us turn now to the analogous relations between the modules. The set 
$A^0_{\tau}$ contains a minimal element, namely $I^0_{\tau}=[0,h_{\tau})$,
which is the smallest interval with $h_{\tau}>0$ and $c_{h_{\tau}}\neq p-1$. 
 We denote the modules
$$
\zer {\cal S}_p^{\tau}\,=\,{\cal M}_{\nu_{I^0_{\tau}}(\tau)}^{\tau}
\qquad\qquad\mbox{and}\qquad\qquad
\ste {\cal S}_p^{\tau}\;=\;  \raise4pt\hbox{${\cal S}_p^{\tau}$}\Big/
\raise-4pt\hbox{$\zer {\cal S}_p^{\tau}$}\;. 
$$
By {\em 3.} of Theorem~\ref{thm-KS99} we know that the modules in the 
composition series of ${\cal M}_{\nu_{I^0_{\tau}}(\tau)}^{\tau}$ are of the
form ${\cal D}_p^{\nu_J(\tau)}$,  where $I^0_{\tau}\subseteq J$, which is equivalent to
$J\in A_{\tau}^0$. Hence 
\beq
{\cal F}(\zer {\cal S}_p^{\tau})\;=\;\{{\cal D}_p^{\nu_I(\tau)}\,:\;I\in A_{\tau}^0\}
\qquad\quad\mbox{and}\qquad\quad
{\cal F}(\ste {\cal S}_p^{\tau})\;=\;\{{\cal D}_p^{\nu_I(\tau)}\,:\;I\in A_{\tau}^+\}\;.
\eeq
With $h_{\tau'}=k_{\tau}$ and $I_{\tau'}^0=\phi(\emptyset)$ we thus have 
$\zer {\cal S}_p^{\tau'}={\cal M}^{\tau'}_{\tau}$. The bijection from 
Corollary~\ref{cor-bij} allows us now to show that 
the composition series of $\zer {\cal S}_p^{\tau'}$ and $\ste {\cal S}_p^{\tau}$ 
yield exactly the same set of irreducibles 
via the identification 
\beq\label{eq-phistar}
\phi^*\;:\;\;{\cal F}(\zer {\cal S}_p^{\tau'})\;\stackrel{=}{\longrightarrow}
\;{\cal F}(\ste {\cal S}_p^{\tau})\;\;:\;\;\;\;{\cal D}^{\nu_I(\tau')}_p
\;\,\stackrel{=}{\mapsto}\,\;{\cal D}^{\nu_{\phi(I)}(\tau)}_p\;. 
\eeq

\begin{lemma}\label{lm-maps} Let $\tau$ and $\tau'$ be as above. 
Suppose there is a non-zero map
$$
\xi\;:\;\;{\cal S}_p^{\tau}\;\longrightarrow\;{\cal S}_p^{\tau'}\qquad\mbox{with}
\qquad\zer {\cal S}_p^{\tau}\subseteq ker(\xi)\;. 
$$ 

Then we have 
$$
im(\xi)=\zer {\cal S}_p^{\tau'}\qquad\qquad\mbox{and}\qquad\qquad
ker(\xi)=\zer {\cal S}_p^{\tau}\;. 
$$
\end{lemma}

{\em Proof:} Consider the map $\overline{\xi}:\ste {\cal S}_p^{\tau}\longrightarrow
{\cal S}_p^{\tau'}$ defined on the quotient by $\zer {\cal S}_p^{\tau}$ as well
as the composite $\overline{\overline{\xi}}:\ste {\cal S}_p^{\tau}\longrightarrow
\ste {\cal S}_p^{\tau'}$ with the projection. Now, since by {\em 1.} of 
Theorem~\ref{thm-KS99} we have a {\em disjoint} union 
${\cal F}({\cal S}_p^{\tau'})={\cal F}(\zer {\cal S}_p^{\tau'})\cup 
{\cal F}(\ste {\cal S}_p^{\tau'})$ so that by (\ref{eq-phistar})
${\cal F}(\ste {\cal S}_p^{\tau'})\cap {\cal F}(\ste {\cal S}_p^{\tau})
=\emptyset$. Hence $\overline{\overline{\xi}}$ is a map between modules
with no common irreducibles in their composition series so that 
$\overline{\overline{\xi}}=0$. Thus the image of $\xi$
or $\overline{\xi}$  must lie  in $\zer {\cal S}_p^{\tau'}$. Consider the
sequence of maps
$$
\xi\;:\;\;\;{\cal M}^{\tau}_{\tau}={\cal S}_p^{\tau}\;
\stackrel{\gamma}{\longrightarrow}\;
\ste {\cal S}_p^{\tau}\;
\stackrel{\overline{\xi}}{\longrightarrow}\;
\zer {\cal S}_p^{\tau'}={\cal M}^{\tau'}_{\tau}\; 
$$
Since  all multiplicities are one we have again a disjoint union 
${\cal F}({\cal M}^{\tau}_{\tau})={\cal F}(ker(\xi))\cup 
{\cal F}(im(\xi))$ so that either ${\cal D}_p^{\tau}\in {\cal F}(ker(\xi))$
or ${\cal D}_p^{\tau}\in {\cal F}(ker(\xi))$. In the first case this 
would mean that $ker(\xi)\subseteq {\cal M}^{\tau}_{\tau}$ is a submodule
that contains ${\cal D}_p^{\tau}$ in its composition series so that by
minimality $ker(\xi)= {\cal M}^{\tau}_{\tau}$. We had, however,
assumed that $\xi\neq 0$. Hence  ${\cal D}_p^{\tau}$ must be contained in
the composition series of $im(\xi)\subseteq {\cal M}^{\tau'}_{\tau}$.
Again it follows by minimality that  
$im(\xi)={\cal M}^{\tau'}_{\tau}$. This also means that 
$\overline{\xi}$ is a surjective map. The fact that all 
multiplicities are one together with (\ref{eq-phistar}) implies that
$dim(\ste {\cal S}_p^{\tau})=dim(\zer {\cal S}_p^{\tau'})$. Hence 
$\overline{\xi}$ must be an isomorphism. 
\ep 
\smallskip

We are finally in the position to establish following resolution of 
irreducible modules by Specht modules. 

\begin{theorem}\label{thm-Snexact}
The sequences (\ref{eq-Cseq}) from Corollary~\ref{cor-Cseq} are exact.  
\end{theorem}

{\em Proof:} We begin by deriving from Theorem~\ref{thm-KS99} that the first term
${\cal S}^{\{n+1-2l\}}_p$ in the sequence is irreducible. Here we have 
$b=l$ and $c=n+1-2l=2ph+2q-k-2l$ so that by (\ref{eq-lastindex})
$$
b=\left\{\begin{array}{ll} q-k\;&\;\mbox{if}\;q\geq k\\
                           q \;&\;\mbox{if}\;q<k\\
          \end{array}\right.
\qquad\mbox{and}\qquad
c_0=\left\{\begin{array}{ll} k\;&\;\mbox{if}\;q\geq k\\
                           p-k \;&\;\mbox{if}\;q<k\\
          \end{array}\right.
$$
Now for a set $I\in A_{\tau}$ with $I\neq\emptyset$ we have $\delta_I\geq p-c_0$. 
We also need $ \delta_I\leq b$ and hence  $p-c_0\leq b$. For $q\geq k$ this condition
reduces to $p\leq q$, which is not possible, and for $q< k$ we find $q\geq k$,
a contradiction as well. Hence $A_{\tau}=\{\emptyset\}$ so that
${\cal S}^{\{n+1-2l\}}_p=\ste {\cal S}^{\{n+1-2l\}}_p={\cal D}^{\{n+1-2l\}}_p$
and $\zer {\cal S}^{\{n+1-2l\}}_p=0$. 

Next, we observe that all maps 
$E^{k_i}:{\cal S}^{\{(i+1)p+k_{i+1}\}}_q\to {\cal S}^{\{ip+k_i\}}_q$ in the 
sequence are of the type of the ones in Lemma~\ref{lm-maps} with diagrams
related by  (\ref{eq-tauprime}). More precisely, we have $c_0=k_{i+1}$ so that
$c'=c-2c_0=(i+1)p+k_{i+1}-2k_{i+1}=jp+(p-k_{i+1})=ip+k_i$. 
Corollary~\ref{cor-Cseq} also implies that these maps are non-zero and
equivariant. It now follows by induction, going from large to small $i$,  that
$$
ker(E^{k_i})\;=\; \zer {\cal S}_q^{\{(i+1)p+k_{i+1}\}}
\qquad\quad\mbox{and}\qquad\quad
im(E^{k_i})\;=\; \zer {\cal S}_q^{\{ip+k_i\}} \;. 
$$
For the first two maps in the sequence this is clear by irreducibility of 
the first module and the fact that the next map is non-zero. 

Once we have proved the relation for $E^{k_{i+1}}$ we know from 
$E^{k_i}E^{k_{i+1}}=0$ that $\zer {\cal S}_q^{\{(i+1)p+k_{i+1}\}}=im(E^{k_{i+1}})
\subseteq ker(E^{k_i})$. We can thus apply Lemma~\ref{lm-maps} to $E^{k_i}$
and infer the statement for $i$ from $i+1$. Hence exactness holds for the
terms before  ${\cal S}^{\{k\}}_p$. 

For this last Specht module we have $c_0=k$ and $c_j=0$ for $j>0$. Thus
$A_{\tau}^+=\{\emptyset\}$ and $A_{\tau}^0=\{[0,j): p^j\leq\frac {n+k+1}2\}$
using that $\delta_{[0,j)}=p^j-k$ and $b=\frac {n-k+1}2$. Particularly, we 
have that $\ste {\cal S}^{\{k\}}_p={\cal D}^{\{k\}}_p$ is already irreducible. 
The kernel of the last map already contains $\zer {\cal S}^{\{k\}}_p$, the 
image of the previous map, in its kernel. If the kernel was bigger the map
would thus have to be zero, which is not the case. Hence exactness also holds
at ${\cal S}^{\{k\}}_p$. Finally, the last map is by irreducibility onto
so that exactness holds through out the sequence. 
\ep

As a result we obtain that the sequence of TQFT's in (\ref{eq-seqTQFT})
of Corollary~\ref{cor-seqTQFT} is exact and this yields a resolutions of
the TQFT's $\dov{\cal V}_{p}^{(k)}$ with $0<k<p$, thus proving 
Theorem~\ref{thm-TQFTreso}. The kernels and cokernels of these
sequences also define TQFT's, which we denote in analogy to the symmetric group
representations by $\zer {\cal V}^{(j)}_p$ and $\ste {\cal V}^{(j)}_p$ so 
that we have short exact sequences
\beq\label{eq-zersteV}
0\;\to\;\zer {\cal V}^{(j)}_p\;\longrightarrow\;
 {\cal V}^{(j)}_p\;\longrightarrow\;
\ste {\cal V}^{(j)}_p\;\to\; 0\;\;. 
\eeq

\bigskip

\head{6. Characters, Dimensions, and the Alexander Polynomial}\label{S6}
\nopagebreak

An obvious application of Theorem~\ref{thm-Snexact} is that we can express the
characters $\varphi^{\tau_k}_p$ and dimensions of the ${\cal D}_p^{\tau_k}=
{\cal D}_p^{\{k\}}$ for diagrams 
$\tau_k=[\frac {n+k-1}2, \frac {n-k+1}2]$ with $k\equiv n+1\mod 2$ and $0<k<p$
in terms of the ordinary characters $\chi^{\tau}$ of Specht modules
${\cal S}^{\tau}$.

\begin{coro}\label{cor-characters} For $k\equiv n+1\mod 2$ and $0<k<p$
we have the following identity of $S_n$-characters 
\beq\label{eq-characters}
\varphi^{\tau_k}_p\;\;=\;\;\sum_{i\geq 0} (-1)^{i} \chi^{\tau_{j_i}}\;, 
\eeq  
where $j_i=ip+k_i$ and $k_i=k$ for $i$ even and $k_i=p-k$ for $i$ odd. 
\end{coro}
As an example we consider the special case $p=5$, where $k=1$ or $3$ if $n$ is
even and $k=2$ or $4$ if $n$ odd. In \cite{Ryba94} Ryba constructs a family of 
irreducible so called
Fibonacci representations $R_n$ and $R_n'$ of $S_n$ over $\F_5$ with 
Brauer characters $\varphi_n$ and $\varphi'_n$ respectively. It follows by 
straight forward computation from (\ref{eq-characters}) and the formulae in
Definition~2 of \cite{Ryba94} that $\phi'_n=\phi^{[r,r]}_5$ and  
 $\phi_n=\phi^{[r+1,r-1]}_5$ if $n=2r$ is even, and  
$\phi'_n=\phi^{[r+2,r-1]}_5$ and  
 $\phi_n=\phi^{[r+1,r]}_5$ if $n=2r+1$ is odd. 

\begin{coro}\label{cor-ryba} 
The Fibonacci representations from \cite{Ryba94} are
$$
R_n\;\cong\;\left\{\begin{array} {cl}
        {\cal D}_5^{[r+1,r-1]}&\;\;\mbox{if}\;\;n=2r\\
        {\cal D}_5^{[r+1,r]}&\;\;\mbox{if}\;\;n=2r+1\\
         \end{array}\right. 
\qquad \mbox{and}\qquad 
R_n'\;\cong\;\left\{\begin{array} {cl}
        {\cal D}_5^{[r,r]}&\;\;\mbox{if}\;\;n=2r\\
        {\cal D}_5^{[r+2,r-1]}&\;\;\mbox{if}\;\;n=2r+1\\
         \end{array}\right. 
$$
\end{coro}
We expect similar relations with the generalizations of these representations 
obtained by Kleshchev in \cite{Kle96}. 

The dimensions of the Specht modules are naturally given by the 
{\em Catalan numbers} $C(n,j)={n\choose j}-{n\choose {j-1}}$. More precisely,
$dim({\cal S}^{[n-b,b]})=C(n,b)$. Particularly, we find for the components of
the sequence in (\ref{eq-Cseq}) that  
$dim({\cal S}^{\{c_{2s}\}})=C(n,b -sp)$ and 
 $dim({\cal S}^{\{c_{2s+1}\}})=-C(n,b +(s+1)p)$, 
with $b=\frac {n+1-k} 2\,$,
where we also 
use that $C(n,j)=-C(n,n+1-j)$. The alternating sum of the Specht module dimensions
comes out to be
\beq\label{eq-DdimCat}
d^n_k\;=\;dim({\cal D}^{[n-b,b]}_p)\;\;=\;\;\sum_{s\in\Z} C(n,b+sp)
\qquad\quad\mbox{provided}\quad 0<k= n-2b+1 < p\;. 
\eeq
Note that if we extend the above formula for $d^n_b$ to the next indices we find
\beq\label{eq-dimbd}
d^n_0=0\;\;\mbox{for odd }n\;,\qquad\mbox{and}\qquad d^n_p=0
\;\;\mbox{for even }n
\eeq
In order to describe generating functions for these dimensions we introduce some 
notation. First we write
 $\displaystyle [n]_x=\frac {x^n-x^{-n}}{x-x^{-1}}\in\Z[x,x^{-1}]$ for the usual
quantum integers. Moreover we consider now the ring of cyclotomic integers 
$\Z[\zeta_p]$ obtained $\Z[x,x^{-1}]$ by imposing  the relation 
$\sum_{j=0}^{p-1}x^j=0$. It is free as a  $\Z$-module of rank $p-1$.
We also denote by 
$\Aa [\zeta_p]\subset \Z[\zeta_p]$ the subring invariant under
conjugation $\zeta_p\mapsto\zeta_p^{-1}$. This is a free $\Z$-module
of rank $\frac {p-1} 2$. It is not hard to see that the set of $[k]_{\zeta_p}$ 
restricted to either even $k$ or to odd $k$ yields a $\Z$-basis for 
$\Aa [\zeta_p]$. 
\begin{lemma} \label{lm-dimgenfct}
We have the following identity in $\Aa[\zeta_p]$:
\beq\label{eq-dimgenfct}
[2]_{\zeta_p}^n\;\;=\;\;\sum_{ \mbox{$
\stackrel{0<k<p}{\scriptstyle k\equiv n+1\mod 2}$}}d^n_k[k]_{\zeta_p}
\eeq
\end{lemma}
{\em Proof:} We have by Schur Weyl duality that
 $L^n_{\Z}\cong V_2^{\otimes n}\cong\bigoplus_{j\equiv n-1\mod 2} V_j\otimes
{\cal S}^{[\frac {n+j-1}2,\frac {n-j+1}2 ]}$. The operator $x^H$ with
$H\in\sl(2,\Z)$ is well defined and has trace $tr_{V_j}(x^H)=[j]_x$. 
It thus follows that $[2]_x^n=\sum_{l\geq 1, l\equiv n+1\mod 2}[l]_xC(n,\frac{n-l+1}2)$.
Now any such $l$ can be uniquely written in the 
form $l=k+2sp$ or $l=-k+2(s+1)p$ with $s\geq 0$ and $k\equiv n+1\mod 2$ and $0<k<p$. 
 Specializing to a root of unity $x=\zeta_p$ we have then that 
$[l]_{\zeta_p}=[k]_{\zeta_p}$ and $[l]_{\zeta_p}=-[k]_{\zeta_p}$ respectively.
Also $C(n,\frac{n-l+1}2)= C(n,\frac{n-k+1}2-sp)$ and 
$C(n,\frac{n-l+1}2)= C(n,\frac{n+k+1}2-(s+1)p)=-C(n,\frac{n-k+1}2+(s+1)p)$ 
respectively. Hence, the terms for a fixed $k$ are given by
$[k]_{\zeta_p} C(n,b-sp)$ and $[k]_{\zeta_p} C(n,b+(s+1)p)$
for $s\geq 0$
and $b=\frac {n+1-k} 2\,$,
which with $s\in\Z$ adds up to the expression in (\ref{eq-DdimCat}).
\ep 

Using $[2]_x[k]_x=[k-1]_x+[k+1]_x$ and the bases of $\Aa[\zeta_p]\subset \Z[\zeta_p]$
we readily derive from  (\ref{eq-dimgenfct}) the recursion
$d^{n+1}_k=d^{n}_{k-1}+d^{n}_{k+1}$. This translates for $0\leq a-b\leq p-2$
to
\beq\label{eq-dimrec}
dim({\cal D}_p^{[a,b]})\;=\;\left\{
\begin{array}{cl}   dim({\cal D}_p^{[a,a-1]})\;&\;\;\mbox{if}\;\;a=b\\
dim({\cal D}_p^{[a-1,b]})+dim({\cal D}_p^{[a,b-1]})\;&\;\;\mbox{if}\;\;0< a-b< p-2\\
dim({\cal D}_p^{[a-1,b]})\;&\;\;\mbox{if}\;\;a-b=p-2\\
\end{array}\right.\;. 
\eeq
It is easy to see from this form that the dimensions are indeed given by the
number of paths through the set of diagrams with $a-b\leq p-2$ as described in 
in \cite{Mat96} for general diagrams. 

In the case $p=5$ this recursion reduces to $dim(R_n)=dim(R_{n-1})+dim(R_{n-1}')$ 
and $dim(R_n')=dim(R_{n-1})$  so that  the dimensions are given by 
 Fibonacci numbers. More precisely, we have $dim(R_n)=f_n$ and 
$dim(R_n')=f_{n-1}$, where $f_n$ are the Fibonacci numbers defined by 
$f_0=0$, $f_1=1$ and $f_{n+1}=f_n+f_{n-1}$. 
Note that together with (\ref{eq-DdimCat}) we find interesting presentations 
of Fibonacci numbers in terms of alternating, 5-periodic  sums of Catalan numbers:
\begin{eqnarray}
f_{2r}&=&C(2r,r-1)-C(2r,r-3)+C(2r,r-6)-C(2r,r-8)+\ldots\label{eq-fibcatideven}\\
&=&C(2r+1,r-1)-C(2r+1,r-2)+C(2r+1,r-6)-C(2r+1,r-7)+\ldots\nonumber \\
\mbox{and}\quad&&\nonumber\\
f_{2r+1}&=&C(2r+1,r)-C(2r+1,r-3)+C(2r+1,r-5)-C(2r+1,r-8)+\ldots\label{eq-fibcatidodd}\\
&=&C(2r+2,r+1)-C(2r+2,r-3)+C(2r+2,r-4)-C(2r+2,r-8)+\ldots\nonumber 
\end{eqnarray}
The reader is invited to check these identities independently via recursion 
relations such as $C(n+1,j)=C(n,j)+C(n,j-1)$ or via the well known generating functions
of Catalan and Fibonacci numbers.

Another useful tool in the determination of dimensions are fusion algebras 
or Verlinde algebras, see \cite{FK93}. For the quantum group $U_q(\sl_2)$
at a $p$-th root of unity the ring of irreducible representations yields 
the fusion algebra $\Phi_p$ generated by the irreducibles 
$\lb 1\rb, \lb 2\rb, \ldots, \lb p-1\rb$.  
Let us denote by $mult(k,R)\in\N\cup\{0\}$ the multiplicity of 
$\lb j\rb$ in an element $R\in\Phi_p$ so that $R=\sum_jmult(k,R)\lb k\rb$.
We have relations 
$\lb 1\rb\circ  \lb k\rb= \lb k\rb$, $\lb p-1\rb\circ\lb k\rb= \lb p-k\rb$
and $\lb 2\rb\circ \lb k\rb= \lb k+1\rb+\lb k-1\rb$ for $1<k<p-1$.
Comparing this to the recursion in (\ref{eq-dimrec}) we find that
\beq\label{eq-dimmult}
dim({\cal D}_p^{[\frac{n+k-1}2,\frac{n-k+1}2]})\;=\;mult(k,\lb 2\rb ^{n})
\eeq
Now, as the ${\cal D}_p^{\tau}$ are isomorphic to the weight spaces 
$\dov {\cal W}^{(k)}_p(\lambda, g)$ with $n=n(\lambda)$ we find 
$
dim(\dov {\cal V}^{(k)}_p(\Sigma_g))=\;\sum_{n=0}^g2^{g-n}{g\choose n}
dim({\cal D}_p^{[\frac{n+k-1}2,\frac{n-k+1}2]})$ so that we find the 
following Verlinde type formula. 
\begin{lemma}\label{lm-dimspaces}
\beq\label{eq-dimspaces}
dim(\dov {\cal V}^{(k)}_p(\Sigma_g))\;=\; mult(k,\ff^{g}_p)
\qquad\quad\mbox{with}\;\;\;\;\ff_p=2\lb 1\rb + \lb 2\rb\;\in\,\Phi_p\;\;.
\eeq
\end{lemma}
Compare this to the  TQFT's ${\cal V}_p^{*RT}$ and  
${\cal V}_p^{RT}$  of Reshetikhin Turaev for the quantum groups $U_q(\sl_2)$
and  $U_q({\mathfrak s}{\mathfrak 0}_3)$ respectively 
at a $p$-th root of unity. (${\cal V}_p^{RT}$ is really a factor TQFT 
obtained from ${\cal V}_p^{*RT}$ by restricting to odd dimensional representations).
With $\FF_p=\sum_j\lb 2j+1\rb^2$ and 
$\FF^*_p=2\FF_p=\sum_k\lb k\rb^2$ we have 
\beq\label{eq-verlinde}
dim({\cal V}_p^{*RT}(\Sigma_g))\;=\;mult(1,{\FF^*_p}^g )\;=\;2^g
mult(1,\FF^g_p )\;=\;2^gdim({\cal V}_p^{RT}(\Sigma_g))
\eeq
In the case $p=5$ these formulae allow us to efficiently 
compute and compare dimensions. 
\begin{lemma} We have for the dimensions 
$D^{(k)}_g=dim(\dov {\cal V}_5^{(k-1)}(\Sigma_g))\,$.
\begin{enumerate}
\item For even $g$: 
$$
\begin{array}{rcl}
D^{(1)}_g&=&\frac 12\Bigl(5^{\frac g2} f_{g-1}+f_{2g+1} \Bigr) \\
D^{(4)}_g&=&\frac 12\Bigl(5^{\frac g2} f_{g-1}-f_{2g+1}        \Bigr)\\
\end{array}
\qquad\qquad
\begin{array}{rcl}
D^{(2)}_g&=&\frac 12\Bigl( 5^{\frac g2} f_{g}+f_{2g}       \Bigr)\\
D^{(3)}_g&=&\frac 12\Bigl(  5^{\frac g2} f_{g}-f_{2g}       \Bigr)\\
\end{array}
$$
\item For odd $g$:
$$
\begin{array}{rcl}
D^{(1)}_g&=&\frac 12\Bigl(5^{\frac {g-1}2}(f_{g-2}+f_{g}) +f_{2g+1}  \Bigr)\\
D^{(4)}_g&=&\frac 12\Bigl(5^{\frac {g-1}2}(f_{g-2}+f_{g}) -f_{2g+1}  \Bigr)\\
\end{array}
\qquad\qquad
\begin{array}{rcl}
D^{(2)}_g&=&\frac 12\Bigl( 5^{\frac {g-1}2}(f_{g-1}+f_{g+1}) +f_{2g}      \Bigr)\\
D^{(3)}_g&=&\frac 12\Bigl( 5^{\frac {g-1}2}(f_{g-1}+f_{g+1}) -f_{2g}            \Bigr)\\
\end{array}
$$
\end{enumerate}
and
\beq\label{eq-dimRT5FN}
dim({\cal V}_5^{RT}(\Sigma_g))\;\;=\;\;D^{(1)}_g\,+\,D^{(4)}_g\;\;\qquad\forall g\geq 0\;\;.
\eeq
 
\end{lemma} 

{\em Proof:} Let us use a more convenient notation $1=\lb 1\rb$, $\rho=\lb 3\rb$,
$\sigma=\lb 4\rb$, and $\sigma\circ \rho=\lb 2\rb$ subject to relations
$\rho\circ \rho=1+\rho$ and $\sigma^2=1$. These relations imply 
$\rho^n=f_{n-1}1+f_n\rho$,  $\ff_5=2+\sigma\circ \rho$ and 
$\FF_5=2+\rho$.   We note now that 
$\FF_5^2=(2+\rho)^2=4+4\rho+\rho^2=5(1+\rho)=5\rho^2$ so that 
$\FF_5^{g}=5^{\frac g2}(f_{g-1}+f_g\rho)$ if $g$ is even. From there 
we compute directly that for odd $g$ we have  
$\FF_p^{g}=5^{\frac {g-1}2}((f_{g-2}+f_g)+(f_{g-1}+f_{g+1})\rho)$.
Consider also $\eta=1-\rho$. We have $\eta^2=1-2\rho+\rho^2=2-\rho=1+\eta$
and thus again $\eta^n=f_{n-1}+f_n\eta$. We find $(2-\rho)^g=\eta^{2g}=
f_{2g-1}+f_{2g}\eta=f_{2g+1}-f_{2g}\rho$. 

Now, we have $(1+\sigma)\circ\ff_5=(1+\sigma)\circ\FF_5$ and thus  
$(1+\sigma)\circ\ff_5^g=(1+\sigma)\circ\FF_5^g$. Similarly, we find 
$(1-\sigma)\circ\ff_5^g=(1+\sigma)\circ(2-\rho)^g$ so that 
$\ff_5^g=\frac 12 (1+\sigma)\circ\FF_5^g+\frac 12(1-\sigma)\circ(2-\rho)^g$.
With the previous results on $\FF_5^g$ and $(2-\rho)^g$ we thus find 
a formula for $\ff_5^g$, which inserted into (\ref{eq-dimspaces}) yields the
asserted formulae. 
\ep 

The formula in (\ref{eq-dimRT5FN}) reflects the fact that 
$\dov {\cal V}_5^{(1)}\oplus \dov {\cal V}_5^{(4)}$ is the sum of the irreducible
constituents of the $\F_5$-reductions of ${\cal V}_5^{RT}$ as TQFT's,
see \cite{KerFib}. For larger primes $p$ it is, however, not
possible that a $\F_p$-reduction ${\cal V}_p^{RT}$ has only the 
${\cal V}_p^{(j)}$ as irreducible components. This is easily see by
looking at the large $g$ asymptotics of the dimension expressions 
in (\ref{eq-dimspaces}) and (\ref{eq-verlinde}). The operations on $\Phi_p$ given
by multiplication by $\FF_g$  or $\ff_g$ are represented by matrices with non-negative 
integer coefficents. Perron-Frobenius theory thus implies that the matrix elements of
of $\FF_g^g$  or $\ff_g^g$, such as those in  (\ref{eq-dimspaces}) and (\ref{eq-verlinde}), 
 grow like $\sim\| \ff_p\|^g$ and $\sim\| \FF_p\|^g$, respectively, where 
$\| \FF_p\|=\frac p{4\sin^2(\frac \pi p)}$ and $\| \ff_p\|=4\cos^2(\frac \pi{2p})$
are the largest eigenvalues of the associate matrices. We thus obtain (\ref{eq-PFasy}). 

Note, that $\| \FF_5\|=\| \ff_5\|$ but that  $\| \FF_p\|>\| \ff_p\|$ if $p>5$.
Thus, a linear relation as Theorem~\ref{thm-KerFib} cannot generalize to $p>5$. 
Instead, we can find polynomials $\,R_p(f)\in\Z[f]\,$ of with degree
$\displaystyle {\rm deg}(R_p)=\frac {p-3}2$ such that 
$\FF_p\;=\;R_p(\ff_p)\,$ by using recursive relations in $\Phi_p$, which
may be identified with a ${\mathbb F}_2$-extension of the real part of $\Z[\zeta_p]$. 
Using modified Tschebycheff polynomials, which we define by the recursion 
$\,P_{j+1}(x)+P_{j-1}\,=\,x\cdot P_j(x)\,$, with $P_0(x)=1$ and $P_1(x)=x$,
the $R_p$ can be written as follows.
\beq
R_p(f)\;\;=\;\;\sum_{j=0}^{\frac {p-3}2}\,n_j\cdot P_j(f-2)\;,
\qquad\quad{\rm with}\qquad n_j\,=\,
\left\{
\begin{array}{cl}
\frac {p-1-j}2&\mbox{for $j$ even}\\
\frac {j+1}2&\mbox{for $j$ odd}
\end{array}
\right.\;.
\eeq
For example, 
$R_5(f)=f$, 
$R_7(f)=2f^2-7f+7$, 
$R_9(f)=2f^3-9f^2+9f+3$,
$R_{11}=3f^4-22f^3+55f^2-55f+22$, and 
$R_{13}(f)=3f^5-26f^4+78f^3-91f^2+26f+13$.
These polynomial appear to play an important r\^ole in representation theoretic
aspects of Conjecture~\ref{conj-sums}. 

\medskip

The analogs of the  character expansions given  in 
Corollary~\ref{cor-characters} and relations as in Lemma~\ref{lm-dimgenfct}
in the  context of the corresponding TQFT's
attain a topological interpretation via the Alexander Polynomial 
$\Delta_{\varphi}(M)\in\Z[x,x^{-1}]$ for a compact, oriented 
3-manifold $M$ with a selected epimorphism
$\varphi:H_1(M,\Z)\to\mkern-15mu\to\Z$. Up to $S$-equivalence the cocycle $\varphi$ defines a
two-sided embedded surface $\Sigma\subset M$ and 
cobordism $C_{\Sigma}:\Sigma\to\Sigma$ obtained by removing a neighborhood
of $\Sigma$ from $M$. The Alexander 
Polynomial is then given (up to a sign, which is determined by the 
additional framing structure on $M$)
by the following identity extracted
in Section~11 of \cite{Ker01}:  
\beq\label{eq-AlexClass}
\Delta_{\varphi}(M)\;\;=\;\;trace(x^{-H}{\cal V}_{\Z}(C_{\Sigma}))\;\;=\;\;
\sum_{j\geq 1}[j]_{-x}trace({\cal V}^{(j)}_{\Z}(C_{\Sigma}))
\eeq
By inserting $x=\zeta_p$ or  $x=-\zeta_p$, a $p$-th root of unity, and 
reducing the integer
coefficients to $\F_p$ we consider the image of the Alexander Polynomial 
under the two so defined natural maps
\beq\label{eq-defredAlex}
\Z[x,x^{-1}]\;\longrightarrow\;\F_p[\zeta_p]\;\;\;:\;\;\;x\;\mapsto\;\pm\zeta_p\;\;:
\Delta_{\varphi}(M)\;\mapsto\;\dov \Delta_{\varphi, p}^{\pm}(M)\;\;. 
\eeq
Now, if we insert $x=\zeta_p$ into (\ref{eq-AlexClass}) and use
$[j_i]_{\pm\zeta_p}=(\pm 1)^{k-1}(-1)^i[k]_{\zeta_p}$, where $j_i=ip+k_i$ as 
before, we infer from the 
resolutions of TQFT's in (\ref{eq-seqTQFT}) the following identities
\begin{eqnarray}\label{eq-alexalt}
\dov \Delta_{\varphi, p}^{+}(M)
\;\;&=&\;\;\sum_{k=1}^{\frac{p-1}2}(-1)^{k-1}[k]_{\zeta_p}\Bigl(
trace(\dov {\cal V}^{(k)}_p(C_{\Sigma}))+
trace(\dov {\cal V}^{(p-k)}_p(C_{\Sigma}))\Bigr)\;, \\
\dov \Delta_{\varphi, p}^{-}(M)
\;\;&=&\;\;\sum_{k=1}^{\frac{p-1}2}\,[k]_{\zeta_p}\Bigl(
trace(\dov {\cal V}^{(k)}_p(C_{\Sigma}))-
trace(\dov {\cal V}^{(p-k)}_p(C_{\Sigma}))\Bigr)\;. 
\end{eqnarray}
This implies Theorem~\ref{thm-pAlexander}. At a 5-th root of unity 
$\dov \Delta_{\varphi, p}^{+}(M)$ depends only on traces of $C_{\Sigma}$ 
under ${\cal V}^{(1)}_5\oplus {\cal V}^{(4)}_5$ and 
${\cal V}^{(2)}_5\oplus {\cal V}^{(3)}_5$. The former is in \cite{KerFib} identified 
with the integral semisimple reduction of the Reshetikhin Turaev theory 
${\cal V}^{RT}_5$ and implies (\ref{eq-Lescop}).
\bigskip

\head{7. Johnson-Morita Extensions}\label{S7}
\nopagebreak

As before we consider for  $H=H_1(\Sigma)$  a standard  basis 
$\{a_1,\ldots,a_g,b_1,\ldots,b_g\}$ that is symplectic with respect to 
the standard skew form $(\,,)$ and orthonormal with respect to the inner
form $\lz\,,\rz$. Let $J\in {\rm Sp}(2g,\Z)$  be the special element defined
by $Ja_i=b_i$ and $Jb_i=-a_i$ so that  $\lz x,y\rz=(x,Jy)$ and
$Jg^{-1}J^{-1}=g^*$. Also denote by $\omega=\sum_ja_j\wedge b_j$
the standard invariant 2-form. 
   
For any $x\in\ext mH$ we can now define a degree-$m$ 
map $\nu(x):\ext *H\to \ext {*+m}H$ by $\nu(x).y=x\wedge y$. 
From this we define another map as $\mu(x)=\nu(Jx)^*:
\ext *H\to \ext {*-m}H$ of degree $-m$. 

\begin{lemma} The maps $\nu$ and $\mu$ have the following properties: 
\begin{enumerate}
\item Covariance: $g\nu(x)g^{-1}=\nu(gx)$
and $g\mu(x)g^{-1}=\mu(gx)$ for all $x\in\ext * H$ and $g\in {\rm Sp}(2g,\Z)$. 
\item Homomorphy: $\nu(x\wedge y)=\nu(x)\nu(y)$  and $\mu(x\wedge y)=\mu(y)\mu(x)$.
\item Generators: $\nu(\omega)=E$ and $\mu(\omega)=F$ so that $[E,\nu(x)]=[F,\mu(x)]=0$ for all $x\in\ext *H$. 
\item Anticommutator: $\mu(x)\nu(y)+\nu(y)\mu(x)=(x,y)\id$ for all $x,y\,\in\,H\,.$. 
\item Commutators: $\displaystyle 
\qquad[E,\mu(x)]=\nu(x)\qquad\mbox{and}\qquad [F,\nu(x)]=-\mu(x)
\qquad\mbox{for}\;\;x\,\in\,H\,.$
\end{enumerate}
\end{lemma}

{\em Proof:} Covariance for $\nu$ is obvious. For $\mu$ consider
$g\mu(x)g^{-1}=(g^{*-1}\nu(Jx)g^*)^*=\nu(g^{*-1}Jx)^*=
\nu(Jgx)^*=\mu(gx)$. Also homomorphy is obvious and the fact that $\nu(\omega)=E$
follows by definition. This implies the zero-commutators since $\omega$ is
even, hence central,  and $F=E^*$. The anticommutator relation is 
readily translated to $\nu(x)^*\nu(y)+\nu(y)\nu(x)^*=\lz x,y\rz\id$.
As a symmetric bilinear relation it suffices to prove this for a system of two
orthonormal vectors $v$ and $w$. Then $\ext *H=\ext *L\otimes \ext *L^{\perp}$,
where $L$ is the space spanned by $v$ and $w$ and $ L^{\perp}$ its orthogonal 
complement. Clearly, $\nu(v)$ and $\nu(w)$ act only on the  
$\ext *L\cong({\Z}^2)^{\otimes 2}$ with
basis $\{ 1,v,w,v\wedge w\}=\{e_-\otimes e_-, e_+\otimes e_-, 
e_-\otimes e_+, e_+\otimes e_+\}\,$. In this form the operators take on the
form $\nu(v)=E\otimes\id$ and $\nu(w)=-H\otimes E$, where the $\sl_2$-generators
$E$ and $H$ act as usual. The relation follows now from $HE+EH=0$ and
$E^*E+EE^*=FE+EF=\id$. The commutators follow by direct calculation. 
$E\mu(x)=\nu(\omega)\mu(x)=\sum_i\nu(a_i)\nu(b_i)\mu(x)=
-\sum_i\nu(a_i)\mu(x)\nu(b_i)+\sum_i\nu(a_i)(x,b_i)=
\sum_i\mu(x)\nu(a_i)\nu(b_i)  -\sum_i(x,a_i)\nu(b_i)  +\sum_i\nu(a_i)(x,b_i)
=\mu(x)E+\nu\Bigl( \sum_ia_i(x,b_i)-b_i(x,a_i)\Bigr)
=\mu(x)E+\nu\Bigl( \sum_ia_i\lz x,a_i\rz+b_i\lz x,b_i\rz\Bigr)
=\mu(x)E+\nu(x)$. The second relation is just the  conjugate of the first. 
\ep 
\smallskip

Equipped with these relations we can now construct maps between the 
${\rm Sp}(2g,\Z)$-representation spaces from previous sections.

\begin{lemma}\label{lm-muall} 
The map $\mu$ restricts and factors into an ${\rm Sp}(2g,\Z)$-covariant map 
$$
\mu\;:\;\;
\frac {\ext m H}{\omega\wedge\ext {m-2} H}
\;\;\longrightarrow\;\;Hom({\cal V}^{(j)}_{\Z}(\Sigma),
{\cal V}^{(j+m)}_{\Z}(\Sigma))\;.  
$$
Moreover, for any $x\in\ext m H$ we have that 
\beq\label{eq-muEmap}
\mu(x): im(E^{l+m})\;\longrightarrow\; im(E^{l})\;. 
\eeq
\end{lemma}

{\em Proof:} Recall that ${\cal V}^{(j)}(\Sigma_g)=ker(F)\cap \ext {(g-j+1)}H$. 
Now $\mu(x)$ commutes  with $F$ and thus maps $ker(F)$ to itself and thus by 
counting degrees to ${\cal V}^{(j+m)}(\Sigma_g)=ker(F)\cap \ext {(g-j-m+1)}H$. 
Moreover, $\mu(\omega\wedge x)=\mu(x)\mu(\omega)=\mu(x)F=0$ if restricted to
$ker(F)$. By the homomorphism property it suffices to show the second relation for
$m=1$. In this case $\mu(x)E^{l+1}z=E^{l+1}\mu(x)z+(l+1)E^l\nu(x)z
\in im(E^l)$ by iteration
of the commutator relation. 
\ep 
\smallskip

Consider now the following family of extensions of the symplectic  group 
 defined as a semidirect product. 
\beq\label{eq-defJMgroup}
JM_{a}(m,g)=\Bigl( \frac 1 a \frac {\scriptstyle \sext m H}
{\,\, \omega\wedge\sext {m-2}\scriptstyle H\,}\Bigr)\rtimes {\rm Sp}(2g,\Z)
\eeq

\begin{propos}\label{propos-Umods} For $a\not\equiv 0\mod p$ 
we have well defined  representations of $JM_{a}(m,g)$ on 
$$
{\cal U}^{(j)}_p(m,g)\;\;=\;\;
{\cal V}_p^{(j)}(\Sigma_g)\stackrel{\mu}{\oplus}{\cal V}_p^{(j+m)}(\Sigma_g)
$$
which decomposes as indicated if restricted to ${\rm Sp}(2g,\Z)$ and is given by $\mu$ if
restricted to the abelian part. 
For $j\equiv k\mod p$ with  $0<k<p-m$ and using the notation from
(\ref{eq-zersteV}) we have that 
$\zer {\cal U}^{(j)}_p(m,g)\;\;=\;\;
\zer {\cal V}_p^{(j)}(\Sigma_g)\stackrel{\mu}{\oplus}
\zer {\cal V}_p^{(j+m)}(\Sigma_g)$ is a proper submodule with subquotient 
\beq\label{eq-Umods}
\ste {\cal U}^{(j)}_p(m,g)\;\;=\;\;
\ste{\cal V}_p^{(j)}(\Sigma_g)\stackrel{\mu}{\oplus}
\ste {\cal V}_p^{(j+m)}(\Sigma_g)\;\;\cong\;\;
\raise3pt\hbox{${\cal U}^{(j)}_p(m,g)$}\Big/ 
\raise-3pt\hbox{$\zer {\cal U}^{(j)}_p(m,g)$}\;. 
\eeq
\end{propos}

{\em Proof:} In general if $V$, $W$ and $M$ are $G$-modules, and 
$\mu:M\to Hom(V,W)$ a covariant
map,   we construct a module $V\stackrel{\mu}{\oplus}W$ of 
 $M\rtimes G$  by letting $(m,g)$ act on the sum $V\oplus W$ by the 
block matrix $\scriptstyle \left[\matrix {g & 0\cr \mu(m)g &g\cr}\right]$
so that $W\subset V\stackrel{\mu}{\oplus}W$ is a submodule with subquotient $V$. 
Thus the map $\mu$ from Lemma~\ref{lm-muall} defines such a module 
for $JM_1(m,g)$ over $\Z$. In the $\F_p$-reduction  $a$ is invertible so
that $\mu$ can be extended to $JM_a(m,g)$.

By exactness of the sequences in (\ref{eq-seqTQFT}) we have that
 $\zer {\cal V}_p^{(j)}(\Sigma_g)={\cal V}_p^{(j)}(\Sigma_g)\cap im(E^{p-k})$.
Thus by (\ref{eq-muEmap}) of Lemma~\ref{lm-muall} we have that 
$\mu(x)(\zer {\cal V}_p^{(j)}(\Sigma_g))\subseteq 
{\cal V}_p^{(j+m)}(\Sigma_g)\cap im(E^{p-k-m})=
\zer {\cal V}_p^{(j+m)}(\Sigma_g)$. This implies, by construction,
that  $\zer {\cal V}_p^{(j)}(\Sigma_g){\oplus}
\zer {\cal V}_p^{(j+m)}(\Sigma_g)$ is indeed a submodule.
\ep
\smallskip

Note that for $0<k<p$ the factors in (\ref{eq-Umods}) are irreducible. 
We also write
\beq\label{eq-Uirred}
\dov {\cal U}^{(j)}_p(m,g)\;\;=\;\;
\dov{\cal V}_p^{(j)}(\Sigma_g)\stackrel{\mu}{\oplus}
\dov {\cal V}_p^{(j+m)}(\Sigma_g)\;. 
\eeq
The case that is topologically relevant is $m=3$. For this case Morita
constructed in \cite{Mor93} a homomorphism  $\tilde k:\Gamma_g\to JM_2(3,g)$
on the mapping class group $\Gamma_g$. Its kernel is the group  ${\cal K}_g$
generated by bounding cycles and its restriction to the Torelli group coincides
with the Johnson homomorphism $\tau_2:{\cal I}_g\to\frac {\sext 3 H}{ H}$ 
from \cite{Joh80}. 

\begin{theorem}[Johnson\cite{Joh80}, Morita \cite{Mor93}] 
There is a finite index subgroup 
${\bf Q}_g\subset JM_2(3,g)=  \frac 1 2 \frac { \sext 3 H}
{\,\,  H\,} \rtimes {\rm Sp}(2g,\Z)$, and  homomorphisms $\tau_2$ and $\tilde k$
such that  the following diagram is commutative with all rows and the last 
column exact. 
\begin{diagram}[height=0.7cm,width=.9cm]
&&&&&&\;\,\quad &0 &&\\
&&&&&&\qquad &\dMyTo &&\\
0 & \rMyTo &{\cal K}_g & \rMyInto
 &{\cal I}_g& \rMyOnto^{\,\, \tau_2\!\! }&
& \,\frac {\sext 3 H}{\omega \wedge H}\, &\rMyTo &0\\
&&   \dEq && \dMyInto & && \dMyInto &&\\
0 & \rMyTo &{\cal K}_g & \rMyInto &\Gamma_g& \rMyOnto^{\,\, \tilde k\!\! }&
&\,\, {\bf Q}_g\,\, &\rMyTo &0 \\
&&   \dMyInto && \dEq && & \dMyOnto \\
0 & \rMyTo &{\cal I}_g & \rMyInto &\Gamma_g& \rMyOnto &
& \,{\rm Sp}(2g,\Z)\, &\rMyTo &0\\
&&&&&&&\dMyTo &&\\
&&&&&&&0 &&\\
\end{diagram}
\end{theorem}
Combining this result with Proposition~\ref{propos-Umods} now yields 
Theorem~\ref{thm-JMext}, where we denoted the special module 
$\dov {\cal U}^{(j)}_p(\Sigma_g)=\dov {\cal U}^{(j)}_p(3,g)$. Note, that
we obtain from the sequences in (\ref{eq-seqTQFT}) similar resolutions.
Particularly, for $g\geq 0$,    $p\geq 5$ a prime  and $0<k<p-3$
there is an exact  sequence of maps of $\F_p$-modules as follows. 
\begin{eqnarray}\label{eq-seqU}
\ldots\,\longrightarrow&&{\cal U}^{((i+1)p+k_{i+1})}_p(\Sigma_g)\,\longrightarrow\,
{\cal U}^{(ip+k_i)}_p(\Sigma_g)\,\longrightarrow\,\ldots\,\longrightarrow\,
{\cal U}^{(2p+k)}_p(\Sigma_g)\,\longrightarrow\nonumber \\
&\,&\longrightarrow\,\,{\cal U}^{(2p-k-3)}_p(\Sigma_g)\,
\longrightarrow\,
{\cal U}^{(k)}_p(\Sigma_g)\,\longrightarrow
\,\dov {\cal U}^{(k)}_p(\Sigma_g)\,\longrightarrow\,0\;, 
\end{eqnarray}
where now $k_i=k$ for $i$ is even, and $k_i=p-i-3$ if $i$ is odd. 

The maps are alternatingly given by $E^k\oplus E^{k+3}$ and
$E^{p-k}\oplus E^{p-k-3}$. Note, however, that the module
extensions work alternatingly in opposite ways so that the
maps cannot be  ${\bf Q}_g$-equivariant. It is true that by setting 
$\tilde R(g)=R(\hat Jg^{-1}\hat J^{-1})^*$, where $\hat J\in\Gamma_g$ 
is a representative of the $J\in {\rm Sp}(2g,\Z)$ above we can reverse the
exact sequences that define an extension. Yet, even by flipping 
every second extension in (\ref{eq-seqU}) still does not yield 
equivariant maps. 
\medskip

\emptystuff{

\head{8. Open Questions}\label{S8}
\nopagebreak

\begin{enumerate}
\item {\em What is the precise submodule structure of the ${\cal V}_p^{(j)}$?} To answer this
question the results from \cite{KS99} should be combined with the $\sp(2g,\Z)$ and
handle mappings between the Specht modules as worked out explicitly in 
Lemma~\ref{lm-spechtTQFTgen} of Section~3. 
\item {\em Are there generalizations to the resolutions in Corollary~\ref{cor-Cseq} to
general Young diagrams?} Schur Weyl duality suggests to use an analogous 
$\sl (n,\Z)$ action, identify Specht modules as lowest weight spaces, and attempt 
to define differential operators between different modules 
from the  Borel algebra generators $E_i\in \sl (n,\Z)$. Such maps and generalized
resolutions should yields many more character formulae. 
\item {\em Is there an independent proof of exactness?} One should note that we
the property of the modular structure from \cite{KS99} that we really 
used was
that the heads of of some representations are irreducible, 
implying maximal elements and cyclicity. It would be intersting to use this more
directly, possibly by employing also 
 the operators  $\frac 1 {k!}E^k$ that are well defined 
for all $k\geq 0$ on the lattices, though they do not map $ker(F)$'s to each other
are likely to  play a r\^ole. 
Note also that similar such resolutions are of intertest in the representation theory
of $p$-adic groups and quantum groups. 
\item {\em Are there TQFT extensions for all of the $\dov {\cal U}^{(j)}_p(\Sigma_g)$ and 
${\cal U}^{(j)}_{\Z}(\Sigma_g)$ modules?} There are two principal approaches of
constructing these. One is to attempt to enlarge the gauge group used by Frohman 
and Nicas from $U(1)$ to a continuous group $G$  for which $G'=[G,G]\neq 0$ but 
$[G',G]=0$ or at least $[G',G']=0$ such as the Heisenberg group or the group of upper 
triangular $3\times 3$-matrices. The representations of the surface groups 
$\pi=\pi_1(\Sigma)$ will factor then through $N=\frac {\pi}{[\pi,\pi']}$, which 
is the natural space for the Morita extension ${\bf Q}_g$ to act on, see \cite{Mor93}. 
The obvious difficulty  is that such groups are non-compact so that intersection
homology is not well defined. The other line of approach is to consider deformations
of the algebra ${\cal N}$ considered in \cite{Ker01} that would lead to non-trivial
contributions of the Torelli groups.
\item {\em What are the constant order contributions of the Reshetikhin Turaev TQFT's for 
other roots of unity?} We noted in Section~4 that the dimensions of these  vector spaces
grow as $\sim \| \FF_p\|^g$ for large $g$, where $\| \FF_p\|=\frac p{4\sin^2(\frac \pi p)}$,
but the dimensions of the modular components only with $\sim \| \ff_p\|^g$, where 
$\| \ff_p\|=4\cos^2(\frac \pi{2p})$. For example, for $p=7$ these are different numbers 
but they are related by $\,\| \FF_7\|=2\|\ff_7\|^2-7\|\ff_7\|+7\,$. A way of accounting
for quadratic terms is to take {\em tensor products} of TQFT's. These will lead to
${\rm Sp}(2g,\Z)$-representations that are no longer of fundamental weight, and also will 
in general have more components in the modular than the in the integral case. Integer
factor can be accounted for by tensoring with TQFT's, ${\bf C}_n$, for the simple
objects of the underlying category are all invertible forming the 
cyclic group $\Z/n\Z$. We thus expect that the constant order of the 
Reshetikhin Turaev TQFT at a 7-th root of unity is described as a sub and/or quotient
TQFT of some ${\bf C}_2\otimes\!\dov{\cal V}_7^{(\#)}\!\otimes\!\dov{\cal V}_7^{(\#)}$, where 
the $\#$ indicate any linear combination of degrees. In general, $\| \FF_p\|$ is a degree
$\frac {p-3}2$ polynomial in $\| \ff_p\|$ with integer coefficients, which gives some 
indication about this generalizes to higher primes $p$.

\item {\em How do the tensor product TQFT's decompose?} This questions is closely related to
the previous one. Since the Frohman Nicas theories yield ${\rm Sp}(2g,\Z)$-representations 
for all fundamental weights we expect the ring of TQFT's generated by them via tensoring,
summing and decomposing to contain {\em all} possible modular irreducible 
${\rm Sp}(2g,\Z)$-representations that come from heighest weight representations. 

\item {\em How are lowest orders of the Reshetikhin Turaev Invariants related to 
Torsion coefficients?} This considers generalizations of the formula in (\ref{eq-AlexRTid}).
Assuming that the constant order RT-TQFT's are obtained from tensor products 
from the $p$-modular Jacobian TQFT's  we expect by taking traces to 
find corresponding  polynomial expressions in the coefficients of
the Alexander Polynomial. 

\item {\em Which modular parts of the Alexander Polynomial are independent of
the choice of a covering, i.e., 1-cocycle, and why?} The observation made in 
 (\ref{eq-AlexRTid}) derived form the independence of the Reshetikhin Turaev
invariant is expected to have generalizations. Since such independence theorems
are formulated without any reference to quantum topology it is a nearby task
to try to find independent proof entirely within the language of
classical algebraic topology.  Interpretations in terms of 
cyclic covering spaces, generalized signature defects, etc. would be quite useful. 

\item {\em Which criteria for 3-manifold groups can be derived from that?} 
As already noted in the introduction the Alexander Polynomial can be computed for
ant group $G$ with 1-cocycle $\chi:G\to\mkern-15mu\to\Z$. Thus independence of 
$\chi$ of the ${\sf x}^2$-coefficient in (\ref{eq-AlexRTid}) can be thought of as a 
condition on $G$ in order for it to be 3-manifold group $G=\pi_1(M)$. Typically,
this criterion allows us to exclude direct sums of groups with sufficiently different Alexander 
Polynomials as 3-manifold groups, but applies also to more subtle examples. 
Details will be worked out in separate work.

\item {\em What is the relation with other TQFT constructions that relate to the
Alexander Polynomial and Torsion?} Among these are the generalizations of 
Frohman and Nicas to $PSU(n)$-theories \cite{FroNic94}, and constructions 
by Donaldson \cite{Don99} that universally reduce  to the Alexander Polynomial,
or polynomial expressions in its coefficients despite rather intricate 
non-abelian Casson-type constructions. This is somewhat parallel to what we
expect for the constant order expansion of the Reshetikhin Turaev 
invariants as outlined above.

\item {\em 
What are the possible ``deformations'' of these modular TQFT's?}
 Essentially what we mean 
is a classifications of  TQFT's over ${\sf R}_{p,l}=\F_p[{\sf x}]/{\sf x}^{l+1}$ for $l<p$ such that 
their reduction under ${\sf R}_{p,l}\to{\sf R}_{p,0}=\F_p$ is one of the given modular TQFT's. 
The Reshetikhin Turaev theory is in this sense an example of a deformation under the reduction
$\Z[\zeta_p]\to \F_p[\zeta_p]\cong {\sf R}_{p,p-2}$ mapping the root of unity to 
$\zeta_p\mapsto 1+{\sf x}$. The rigidity of the ${\rm Sp}(2g,\Z)$-representation theory and its 
actions on the lower parts of the filtrations of the mapping class groups appear to impose
string constraints on the space of possible deformations.  

\item {\em What are connections to the universal finite type TQFT's? }
The TQFT's discussed here represents the Torelli group non-trivially and thus
have some overlap with the lowest orders of the TQFT extensions of \cite{MO97} to homologically
trivial cobordisms
the famous LMO invariant from \cite{LMO98}.
Similarly, we can consider deformations of the modular TQFT's as described previously and
compare them up to higher  orders of the finite type theory. 
These connections should shed some light on possible modular reductions of the
finite type theories that have a chance to be extended to a homologically 
non-trivial theory.

\end{enumerate} 

}

\bigskip

{\sc\small  The Ohio State University, 

Department of Mathematics,

231 West 18th Avenue,

         Columbus, OH 43210, U.S.A. }

 {\em E-mail: }{ \tt kerler@math.ohio-state.edu}
\end{document}